\documentclass{amsart}

\usepackage{mathrsfs}
\usepackage{multicol}
\usepackage{pdfpages}
\usepackage{wrapfig}
\usepackage{color}
\usepackage{bbm}
\usepackage{tabularx}
\usepackage{caption}
\usepackage{bm}
\usepackage{mathrsfs}
\usepackage{amssymb}
\usepackage{enumitem}
\usepackage{pdfpages}
\usepackage{float}
\usepackage{graphicx}
\usepackage{csquotes}
\usepackage{hyperref}
\usepackage{amsmath}
\usepackage[margin=1.2in]{geometry}
\usepackage{ulem}

\newtheorem{theo}{Theorem}

\newtheorem{lem}{Lemma}

\usepackage{fancyhdr}
\usepackage{ulem}
\usepackage{mathrsfs}
\usepackage[english]{babel}
\usepackage{multicol}
\usepackage{pdfpages}
\usepackage{wrapfig}
\usepackage{color}
\usepackage{tabularx}
\usepackage{caption}
\usepackage{bm}
\usepackage{mathrsfs}
\usepackage{amssymb}
\usepackage{enumitem}
\usepackage{pdfpages}
\usepackage{float}
\usepackage{graphicx}
\usepackage{csquotes}
\usepackage{hyperref}
\usepackage{amsmath}
\usepackage{ulem}
\usepackage{lipsum}

\newtheorem{theorem}{Theorem}[section]

\newtheorem{proposition}[theorem]{Proposition}

\theoremstyle{definition}

\theoremstyle{remark}
\newtheorem{remark}[theorem]{Remark}

\numberwithin{equation}{section}

\numberwithin{equation}{section}
\numberwithin{theorem}{section}
\numberwithin{cor}{section}
\numberwithin{cor}{section}
\numberwithin{eg}{section}
\numberwithin{examp}{section}

\newcommand{\A}{{\mathcal A}}
\newcommand{\bbr}{{\mathcal R}}
\newcommand{\bbn}{{\mathbb N}}
\newcommand{\raw}{\rightarrow}

\newcommand{\nti}{n\rightarrow\infty}

\begin{document}


\title{Central Limit Theorem and Near classical Berry-Esse\'en rate for self normalized sums in high dimensions}

\author{Debraj Das}
\address{Department of Mathematics, Indian Institute of Technology, Bombay, Maharashtra 400076, India}
\email{debrajdas@math.iitb.ac.in}
\thanks{The author was supported in part by DST fellowship DST/INSPIRE/04/2018/001290}






\begin{abstract}
In this article, we are interested in the high dimensional normal approximation of $T_n =\Big(\sum_{i=1}^{n}X_{i1}/\Big(\sqrt{\sum_{i=1}^{n}X_{i1}^2}\Big),\dots,$ $\sum_{i=1}^{n}X_{ip}/\Big(\sqrt{\sum_{i=1}^{n}X_{ip}^2}\Big)\Big)$ in $\bbr^p$ uniformly over the class of hyper-rectangles $\mathcal{A}^{re}=\{\prod_{j=1}^{p}[a_j,b_j]\cap\mathcal{R}:-\infty\leq a_j\leq b_j
\leq \infty, j=1,\ldots,p\}$, where $X_1,\dots,X_n$ are non-degenerate independent $p-$dimensional random vectors. We assume that the components of $X_i$ are independent and identically distributed (iid) and investigate the optimal  cut-off rate of $\log p$ in the uniform central limit theorem (UCLT) for $T_n$ over $\mathcal{A}^{re}$. The aim is to reduce the exponential moment conditions, generally assumed for exponential growth of the dimension with respect to the sample size in high dimensional CLT, to some polynomial moment conditions.    Indeed, we establish that only the existence of some polynomial moment of order $\in [2,4]$ is sufficient for exponential growth of $p$. However the rate of growth of $\log p$ can not further be improved from $o\big(n^{1/2}\big)$ as a power of $n$ even if $X_{ij}$'s are iid across $(i,j)$ and $X_{11}$ is bounded. We also establish near$-n^{-\kappa/2}$ Berry-Esseen rate for $T_n$ in high dimension under the existence of $(2+\kappa)$th absolute moments of $X_{ij}$ for $0< \kappa \leq 1$. When $\kappa =1$, the obtained Berry-Esseen rate is also shown to be optimal. As an application, we find respective versions for component-wise Student's t-statistic, which may be useful in high dimensional statistical inference.
\end{abstract}

\maketitle

\section{Introduction}\label{sec:intro}

Let $X_1,\dots,X_n$ be  independent 
random vectors in $\mathcal{R}^p$ each with mean $0$.  Suppose that $X_i=(X_{i1},\dots$ $,X_{ip})$, $n\geq 1$, and for any $j\in \{1,\dots,p\}$,
 \begin{align*}
S_{nj}=X_{1j}+\dots+X_{nj},\;\;\;\; B_{nj}^2=\sum_{i=1}^{n}\mathbf{E}X_{ij}^2,\;\;\;\; V_{nj}^2 = \sum_{i=1}^{n}X_{ij}^2.   
 \end{align*}
Then the component-wise normalized and the component-wise self-normalized sum of the random vectors $X_1,\dots,X_n$ are respectively defined as
$$U_n = \Big(\frac{S_{n1}}{B_{n1}},\dots,\frac{S_{np}}{B_{np}}\Big)^\prime\;\;\;\; \text{and}\;\;\;\; T_n = \Big(\frac{S_{n1}}{V_{n1}},\dots,\frac{S_{np}}{V_{np}}\Big)^\prime.$$ Define the Gaussian random vector $Z_n$ on $\mathcal{R}^p$ with mean $0$ and covariance matrix $Corr(S_n)$, the correlation matrix of $S_n=\big(S_{n1},\dots,S_{np}\big)^\prime$. Distributional approximation of $U_n$ using the Gaussian random vector $Z_n$  is extensively studied in fixed $p$ setting as well as when $p$ increases with $n$. In this paper we study the Gaussian approximation of $T_n$ when $p$ can grow exponentially with $n$. More precisely, we investigate the rate of growth of $p$ required to ensure the uniform central limit [hereafter referred to as UCLT] given by
\begin{align}
\tau_{n,\A}\equiv   \sup_{A\in \A} \Big| \mathbf{P}(T_n\in A) - \mathbf{P}(Z_n\in A)   \Big|  \raw 0  ~\mbox{as}~ \nti,
    \label{clt-1}
\end{align}
where $\A$ is a suitable collection of convex subsets  
of $\bbr^p$. Typical choices of $\A$  include 

\begin{itemize}
\item
$\A^{dist}
= \Big\{ (-\infty, a_1]\times\ldots\times (-\infty, a_p]: a_1,\ldots, a_p\in \bbr\Big\}$,
the collection of all left-infinite rectangles, leading to the Kolmogorov distance between the distributions of $T_n$ and $Z_n$, 
\item 
$\A^{max} =\Big\{
 (-\infty, t]\times\ldots\times (-\infty, t]: t \in \bbr
\Big\}$, the collection of left infinite hyper-cubes, leading to the Kolmogorov distance between the distributions of $\max_{1\leq j\leq p}T_{nj}$ and $\max_{1\leq j\leq p}Z_{nj}$, 
~~and
\item
$\mathcal{A}^{re}=\Big\{\prod_{j=1}^{p}[a_j,b_j]\cap\mathcal{R}:-\infty\leq a_j\leq  b_j\leq \infty
~\mbox{for}~ j=1,\ldots,p \Big\}$,\\
the collection of 
all hyper rectangles, 
\end{itemize}
among others. Since $\A^{max}\subset \A^{dist}\subset \A^{re}$, we only study the growth rate of $p$ in $\tau_{n,\A^{re}}\rightarrow 0$ as $n\rightarrow \infty$. As an application of the UCLT  (\ref{clt-1}) for the self-normalized random vector $T_n$, we also investigate the UCLT of the component-wise studentized random vector $W_n=(W_{n1},\dots,W_{np})$ over the class of sets $\A^{re}$, where $W_{nj}$ is defined as 
$$W_{nj}=\frac{\sqrt{n}\bar{X}_{nj}}{\sqrt{(n-1)^{-1}\sum_{i=1}^{n}(X_{ij}-\bar{X}_{nj})^2}},$$
with $\bar{X}_{nj}=n^{-1}\sum_{i=1}^{n}X_{ij}$, $j\in \{1,\dots,p\}$. We show that the results that hold for $T_n$ are also true for $W_n$.

Similar to (\ref{clt-1}), we say that UCLT of the normalized random vector $U_n$ holds over a class of sets $\A$ provided 
\begin{align}\label{clt-2}
   \rho_{n,\A}\equiv   \sup_{A\in \A} \Big| \mathbf{P}(U_n\in A) - \mathbf{P}(Z_n\in A)   \Big| \rightarrow 0,\; \text{as}\; n\rightarrow \infty. 
\end{align}
When $p$ is fixed, convergence of $\rho_{n,\A^{re}}$ to $0$ as $n\rightarrow \infty$ follows easily from the classical Lindeberg's  Central Limit Theorem. When $p$ grows with $n$, then there is a series of interesting results available in the literature. In a seminal paper, \cite{CCK13} showed that when the random variables $X_{ij}$'s are sub-exponential, $\rho_{n,\A^{max}}\rightarrow 0$ as $n\rightarrow \infty$ provided $\log p =o(n^{1/7})$. \cite{CCK17} improved their results and established the same growth rate of $p$ when $\A=\A^{re}$ in (\ref{clt-2}).  Under the same assumption of sub-exponentiality of $X_{ij}$'s when $\A=\A^{re}$, the rate of $p$ is improved to $\log p = o(n^{1/5})$ by \cite{CCKK} and \cite{K}. Under the additional assumption that the random vectors $X_i$'s are iid having log-concave density, \cite{FK} improved the growth rate of $p$ to $\log p =o
(n^{1/3}(\log n)^{-2/3})$ in (\ref{clt-2}) with $\mathcal{A}=\mathcal{A}^{re}$.  Recently, \cite{DL} further improved the rate of $p$ to $\log p =o(n^{1/2})$ in (\ref{clt-2}) with $\A=\A^{re}$ when the random vectors $X_i$'s are independent and they have independent \& identically distributed (iid) sub-Gaussian components symmetric around $0$. Some negative results are also available in the literature. For example, \cite{FK} showed that when $X_{ij}$'s are iid across $(i,j)$ with $E(X_{11}^3)\neq 0$ and some other conditions are true, the rate of $p$ in (\ref{clt-2}) with $\A=\A^{max}$ can not be improved further from $\log p =o(n^{1/3})$. On the other hand, \cite{DL} showed using an example that the rate $\log p =o(n^{1/2})$ in (\ref{clt-2}) with $\A=\A^{max}$ can not be improved as a power of $n$ when the underlying setup is symmetric around $0$. Beside finding the optimal growth rate of the dimension $p$ in the UCLT of normalized sum $U_n$, results in the direction of matching the $n^{-1/2}-$rate of  the classical Berry-Esseen theorem have also been investigated in the literature. \cite{FK} established the high dimensional Berry-Esseen rate $(\log p)^{3/2}(\log n)n^{-1/2}$ for $\rho_{n, \mathcal{A}^{re}}$ when $X_i$'s are iid random vectors having log-concave densities. When $X_i$'s are iid random vectors with sub-Gaussian components, \cite{L} established the Berry-Esseen rate $(\log pn)^4 (\log n)n^{-1/2}$ for $\rho_{n, \mathcal{A}^{dist}}$. \cite{KR} generalized the Berry-Esseen rate of \cite{L} to the non-identical case. \cite{FK} also showed that $(\log p)^{3/2}n^{-1/2}$ is the optimal Berry-Esseen rate for $\rho_{n, \mathcal{A}^{re}}$ when $X_{ij}$'s are iid across $(i,j)$ and $EX_{11}^3\neq 0$.


In this paper we study the UCLT as well as establish the near classical Berry Esseen rate for $T_n$ and $W_n$ over the class $\A^{re}$. Here $X_1,\dots,X_n$ are mean zero independent $p-$dimensional random vectors  with each having iid components.  The aim for considering $T_n$ and $W_n$ as an alternative to $U_n$ is to reduce the underlying moment assumptions from existence of some exponential moments to the existence of some polynomial moments of order $\leq 4$. To the best of our knowledge, this is the first work where UCLT is established with exponential growth of the dimension $p$ under only polynomial moment conditions. 
Indeed we show that when $\mathbf{E}|X_{ij}|^{2+\delta}<\infty$ for all $(i,j)$ for some $0<\delta\leq 1$, then the rate $\log p =o\big(n^{\delta/(2+\delta)}\big)$ can be achieved in the UCLT (\ref{clt-1}) of $T_n$ with $\A=\mathcal{A}^{re}$. 
Therefore the rate $\log p = o\big(n^{1/3}\big)$ can be achieved in (\ref{clt-1}) under the existence of third absolute moments of underlying random variables, whereas the underlying random variables are assumed to be sub-exponential to achieve the similar rate of $\log p$ in (\ref{clt-2}) [cf. Corollary 1.1 and Proposition 1.1 in \cite{FK}]. Moreover, the rate of $\log p$ in (\ref{clt-1}) can be improved to $o\big(n^{1/2}\big)$ under the existence of fourth moment and when the third moments vanish. 
Using two examples it is also established that the rates $o\big(n^{1/3}\big)$ and $o\big(n^{1/2}\big)$ are generally optimal for $\log p$ respectively when $X_{ij}$'s are iid having asymmetric distribution and when $X_{ij}$'s are iid having distribution symmetric around $0$. 

When $X_{ij}$'s are all iid, then even no moment condition is required to achieve sub-exponential growth rate of $p$ in (\ref{clt-1}), only the assumption of $X_{11}$ belonging to the domain of attraction of the normal distribution is sufficient. Moreover, all the above mentioned UCLT results still hold if we define $\tau_{n,\A}$ with $T_n$ replaced by $W_n$. The backbone behind establishing the UCLT of $W_n$ from that of $T_n$ is the well-known identity between $T_{nj}$ and $W_{nj}$, given by 
\begin{align}\label{eqn:relation}
\mathbf{P}\Big(W_{nj}\geq x\Big)=\mathbf{P}\Big(T_{nj}\geq x\Big(\frac{n}{n+x^2-1}\Big)^{1/2}\Big),\;\text{for any}\;x\geq 0,  
\end{align}
for all $1\leq j\leq p$ [cf. \cite{E}]. Therefore in the UCLT, similar growth rates of $p$ can be achieved under much weaker moment conditions if we perform component-wise studentization instead of component-wise standardization. This fact is very interesting from the perspective of statistical inference as well. In most of practical problems, the underlying population variance is not known and hence it is essential to perform studentization. 

After establishing UCLTs for self-normalized and studentized sums under different moment conditions, we move towards obtaining different high dimensional Berry-Esseen rates. The main aim is to establish near classical rate for both $T_n$ and $W_n$ over the class $\mathcal{A}^{re}$. We have considered following two cases separately: \\
Case I: $X_{ij}$'s are all iid and $X_{11}$ is in the domain of attraction of the normal distribution.\\
Case II: $X_i =(X_{i1},\dots, X_{ip})^\prime$ are independent random vectors for $i\in \{1,\dots, n\}$. For each $i$, $X_{ij}$'s are iid, $EX_{i1}=0$ and $E|X_{i1}|^{2+\kappa} < \infty$ for some $0< \kappa \leq 1$.\\
We show that for Case I, the high dimensional Berry-Esseen rate is close to the classical Berry-Esseen rate $\delta_{n, 0}$ obtained by \cite{BG}. On the other hand for Case II, we established the high dimensional Berry-Esseen rate which nearly recovers the classical Berry-Esseen rate $d_{n, \kappa}^{-(2+\kappa)}$ obtained by \cite{S05}. The quantities $\delta_{n, 0}$ and $d_{n, \kappa}$ are defined in the beginning of Section \ref{sec:main}. In particular, when $n^{-1}\sum_{i=1}^{n}\mathbf{E}|X_{i1}|^3 = O(1)$ and $\liminf_{n\rightarrow \infty}\big(n^{-1}\sum_{i=1}^{n} EX_{i1}^2\big) >0$ then under the setup of Case II, we show that $(\log(pn))^{3/2}n^{-1/2}$ is the high dimensional Berry-Esseen rate.  We also show that when $\kappa =1$, $(\log p)^{3/2}n^{-1/2}$ is generally the optimal Berry-Esseen rate for self-normalized and studentized random vectors under the setup of Case II. This shows that the obtained Berry-Esseen rate under Case II is generally optimal in high dimensional regimes when $\kappa = 1$.

The proofs of the results crucially  depend on
 a set of non-uniform Berry-Esseen type bounds in the one dimensional CLT for self-normalized sum of  independent random variables
 which, in turn,  heavily depend on the Cram\'er type large deviation for one dimensional self normalized sums. Clearly under the setup of $X_1,\dots, X_n$ being independent with each having iid components,  $X_1,\ldots, X_n$ can also be non-identically 
distributed (e.g. with a different component-wise variance). The independence of $X_1,\dots,X_n$ and the iid nature among the components of each $X_i$ ensure that $T_{n, 1},\dots,T_{n, p}$ are iid. Independence is essential to factorize $\mathbf{P} (T_n \in A)$ over $A \in A^{re}$ into $p$ factors each of which depends only on one dimensional self-normalized sums. Whereas the iid nature of $T_{n, i}$'s is required for the application of Lemma \ref{lem:7}. Lemma \ref{lem:7} simply tells us that it is possible to sort the boundaries of the intervals present in these $p$ factors. This sorting is quite essential in order to define a specific partition of the real line and then to apply a set of suitable error estimates over each range using non-uniform Berry-Esseen bounds for one dimensional self-normalized sum (cf. Lemmas \ref{lem:3}, \ref{lem:5} and \ref{lem:4}). The proof strategy is quite general under the setting considered, in the sense that it is simply reducing the original problem to mere applications of exponential concentration of one dimensional self-normalized sums. The technique is also flexible enough to produce optimal results under different moment conditions (cf. Propositions \ref{prop:2.1}, \ref{prop:2.2} and \ref{prop:4.1}). 

 The literature on one dimensional self-normalized sums is well developed. Preliminary results in self-normalized sums are due to \cite{D}, \cite{E} and \cite{LMRS}. \cite{E} pointed out the crucial relation (\ref{eqn:relation}) which essentially implies that the limit distribution of $T_{n}$ and $W_n$ coincide when $p$ is fixed. Under the assumption that $p=1$ and $X_i$'s are iid, \cite{LMRS} and later \cite{CG04b} found the asymptotic distribution of $T_{n}$ when $X_1$ is in the domain of attraction of some stable law. Uniform Berry-Esseen bounds  for $T_n$ and $W_n$ in one dimension were obtained by \cite{S} and \cite{H} and later refined by \cite{BG} and \cite{BBG}. \cite{WJ} essentially improved these uniform results and established a non-uniform Berry-Esseen bound for $T_n$ and $W_n$. \cite{S97} and \cite{S99} established a Cram\'er type large deviation result for $T_n$ and $W_n$ when $p=1$, which was later improved and generalized by \cite{JSW}. 
 \cite{RW}, \cite{W05}, \cite{WH}, \cite{W11}, \cite{SG} essentially generalized and improved the large deviation results of \cite{JSW}. 
 For an elaborate and systematic presentation of the results for self-normalized sums, one can look into \cite{PLS} and \cite{WS}.

The rest of the paper is organized as follows. We state different UCLT results for self-normalized sums in Section \ref{sec:main}. Respective results corresponding to the component-wise student's t-statistic are presented in Section \ref{sec:student's t}. Section \ref{sec:BE} presents high dimensional Berry-Esseen theorems for both self-normalized and studentized sums.
Proofs of all the results except Proposition \ref{prop:4.1} are presented in Section \ref{sec:pf}. Proofs of Lemma \ref{lem:5}, Lemma \ref{lem:7} and Proposition \ref{prop:4.1} are presented in the Supplementary material file.

\section{UCLT for self-normalized sums}
\label{sec:main}
\setcounter{equation}{0} 
In this section we are going to present different UCLT results for the self-normalized sums in increasing dimension. For the rest of the paper, we assume $\{X_1, \dots, X_n\}$ to be a collection of independent mean zero random vectors in $\mathcal{R}^p$. For each $i\in \{1,\dots,n\}$, $X_i=(X_{i1},\dots,X_{ip})^\prime$ where $\{X_{i1}, \dots, X_{ip}\}$ are iid. We are going to explore the rate of growth of the dimension $p$ for the UCLT  (\ref{clt-1}) over $\A^{re}$ to hold, under different moment conditions. Note that in our setup $Z_n=Z$ for all $n\geq 1$, where $Z$ is the standard Gaussian random vector on $\mathcal{R}^p$ with mean $0$ and covariance matrix ${\mathbb I}_p$, the identity matrix of order $p$. 

We divide this section in two sub-sections. In the first sub-section, we start with exploring the growth rate of $p$ in the general setting when $\mathbf{E}|X_{i1}|^{(2+\delta)}$ are finite for some $0<\delta\leq 1$. Then we assume that $X_{ij}$'s are iid across $(i,j)$ and drop the assumption of finiteness of $(2+\delta)$th absolute moment of $X_{ij}$'s. We only assume that $X_{11}$ belongs to the domain of attraction of the normal distribution. In the second sub-section we explore the rate of $p$ when $\max_{1\leq i \leq n}\mathbf{E}|X_{i1}|^4 = O(1)$, as $n\rightarrow \infty$. 
Before moving to the sub-sections, we are going to define few notations. For a collection of random variables $\{Y_1,\dots, Y_n\}$, define $\sigma_n^2(Y)=n^{-1}\sum_{i=1}^{n}Var(Y_i)$. Also for any $k> 0$, define $\beta_{n,k}(Y)=n^{-1}\sum_{i=1}^{n}\mathbf{E}|Y_i|^{2+k}$ and the Lyapounov's ratio $$d_{n,k}(Y)=\dfrac{\sqrt{n}\sigma_n(Y)}{\big(n\beta_{n,k}(Y)\big)^{1/{(2+k)}}}.$$   
The rate of $p$ is stated in terms of Lyapunov's ratio in Theorem \ref{theo:1} and Theorem \ref{theo:8}. For Theorem \ref{theo:3} and Theorem \ref{theo:7}, define
\begin{align*}
\delta_{n,x}(Y) = \;&n\mathbf{P}\Big(|Y|>\kappa_{n,x}(Y)\Big)+n\big[\kappa_{n,x}(Y)\big]^{-1}\Big|\mathbf{E}\big\{YI\big(|Y|\leq \kappa_{n,x}(Y)\big)\big\}\Big|\\
&+n\big[\kappa_{n,x}(Y)\big]^{-3}\mathbf{E}\big\{|Y|^3I\big(|Y|\leq \kappa_{n,x}(Y)\big)\big\},
\end{align*}
where $\kappa_{n,x}(Y)=\sup\Big\{s:ns^{-2}\mathbf{E}\big\{Y^2I(|Y|\leq s)\big\}\geq 1+x^2\Big\}$, for some random variable $Y$. Whenever $X_{ij}$'s are the underlying random variables, then we simply write $\sigma_n^2$, $\beta_{n,k}$, $d_{n,k}$ and $\kappa_{n,x}$.


\subsection{Rate of growth of $p$ under $(2+\delta)$th absolute moment}

In this sub-section, we state three theorems on the growth rate of $\log p$ in the UCLT over $\mathcal{A}^{re}$. The first two are under the assumption of existence of $(2+\delta)$th absolute moments of $X_{ij}$'s. On the other hand, Theorem \ref{theo:3} shows that sub-exponential growth rate of $p$ is possible in iid case even when we only assume that $X_{11}$ is in the domain of attraction of the normal distribution. 

\begin{theo}\label{theo:1}
Let $\mathbf{E}|X_{n1}|^{2+\delta}<\infty$ for all $n\geq 1$ and for some $0<\delta\leq 1$, $d_{n,\delta}\rightarrow \infty$ as $n\rightarrow \infty$. Then we have $$\tau_{n,\mathcal{A}^{re}}\rightarrow 0 \;\;\; \text{as}\;\;n\rightarrow \infty,$$
provided $\log p =o(d_{n,\delta}^2)$. Moreover, there exists a positive constant $c\leq 1/8$ such that whenever $\log p \leq \epsilon d_{n,\delta}^2$ with $0<\epsilon < c$, $$\limsup_{n \rightarrow \infty}\tau_{n,\mathcal{A}^{re}}< \epsilon^{(1+\delta)/3}.$$
\end{theo}

Theorem \ref{theo:1} shows that
the UCLT \eqref{clt-1} holds with $\A=\A^{re}$
for $p$ growing at the rate $\exp\Big(o\big(d_{n,\delta}^2\big)\Big)$
with the sample size $n$. Moreover, the UCLT holds over $\mathcal{A}^{re}$ for the self-normalized sums with an error even when $\log p$ is exactly  of order $d_{n,\delta}^2$.  In particular, for some $0<\delta\leq 1$ if $m\leq \mathbf{E}|X_{n1}|^{2}, \mathbf{E}|X_{n1}|^{2+\delta}\leq M$ for all $n\geq 1$ for some constants $m,M>0$, then $\log p$ can grow like $o\Big(n^{\delta/{(2+\delta)}}\Big)$ in the UCLT over $\A^{re}$. The rate of growth of $\log p$ obtained in Theorem \ref{theo:1} is in general optimal with respect to the underlying moment conditions, as is shown in the following proposition for $\delta=1$. 
\begin{proposition}\label{prop:2.1}
Let $X_1, \dots, X_n$ be random vectors in $\mathcal{R}^p$ where $X_{i} = (X_{i1},\dots,$ $ X_{ip})^\prime$. Define $\{\beta_n\}_{n\geq 1}$ to be a sequence of real numbers such that $\beta_n \geq e$ and $\beta_n/\sqrt{n}\rightarrow 0$ as $n\rightarrow \infty$. Let $X_{ij}$'s be iid mean zero random variables across $(i,j)$ such that
$$P(X_{11}=-a)=q,\;\;\;\; P(X_{11}=b)= 1- (2\beta_n)^{-1},\;\;\;\;\text{and}\;\; P(X_{11}=2a) = p$$ where $a = \dfrac{b_1}{4} + \dfrac{1}{2}\sqrt{\dfrac{b_1^2}{4}+4\beta_n^2\Big(1-b^2\Big(1-\dfrac{1}{2\beta_n^2}\Big)\Big)}$, $p = \dfrac{1}{6\beta_n^2}\Big(1-\dfrac{b_1}{a}\Big)$, $q = \dfrac{1}{6\beta_n^2}\Big(2+\dfrac{b_1}{a}\Big)$, $b_1=\Big(2\beta_n^2-1\Big)b$ and $b>0$ is sufficiently small. Then if $\dfrac{\log p}{d_{n, 1}^2} \rightarrow \infty$ as $n\rightarrow \infty$, we have $\liminf_{n\rightarrow \infty}\tau_{n, \mathcal{A}^{\max}}>0.$
\end{proposition}

Now there are two directions of improving Theorem \ref{theo:1}. One direction is to establish sub-exponential rate of $p$ even when $(2+\delta)$th absolute  moment does not exist, whereas the other one is to improve the growth rate of $p$ under higher moment conditions. The later one is relegated to the next sub-section. The first direction is considered in the next theorem when $X_{ij}$'s are iid and no moment condition on $X_{11}$ is assumed. The only assumption is that $X_{11}$ is in the domain of attraction of the normal distribution.  

\begin{theo}\label{theo:3}
Let $X_{ij}$'s are iid copies of $X_{11}$ for all $i\in\{1,\dots,n\}$ and $j\in \{1,\dots,p\}$. Let $X_{11}$ be in the domain of attraction of the normal distribution with the distribution of $X_{11}$ being non-degenerate and symmetric around $0$. Define $\omega_n = \big[\delta_{n,0}]^{-1/6}$. 
Then we have
\begin{align*}
\tau_{n,\mathcal{A}^{re}}\rightarrow 0 \;\;\; \text{as}\;\;n\rightarrow \infty, 
\end{align*}
provided $\log p = o(\omega_n^2)$. Moreover there exists a positive constant $c\leq 1/8$ such that whenever  $\log p \leq \epsilon \omega_n^2$ with $0<\epsilon < c$, $$\limsup_{n \rightarrow \infty}\tau_{n,\mathcal{A}^{re}}< \epsilon^{5/3}.$$
\end{theo}
 The above result is remarkable in the sense that sub-exponential growth of the dimension is possible in the UCLT of self-normalized sums even when no moment higher than two exists. Note that $X_{11}$ is in the domain of attraction of normal distribution if and only if $\lim_{t\rightarrow \infty}\big[t^2\mathbf{P}\big(|X_{11}|> t\big)\big]/\big[\mathbf{E}X_{11}^2I\big(|X_{11}|\leq t\big)\big]$ $=0$ (cf. Theorem 4 at page 323 of \cite{CT}). Since due to Lemma 1.3 of \cite{BG}, $\kappa^2_{n,0}/n= \mathbf{E}\big\{X_{11}^2I(|X_{11}|\leq \kappa_{n,0})\big\}$, the necessary and sufficient condition in turn implies that $\delta_{n,0}\rightarrow 0$ as $n\rightarrow \infty$ due to Lemma 1 of C\"org\H{o} et al. (2003). However, $\mathbf{E}|X_{11}|^{2+\delta}$ may not exist for any $\delta>0$ and hence we can not use Theorem \ref{theo:1}. In particular, if we only know that $\mathbf{E}X_{11}^2[\log(1+|X_{11}|)]^{\alpha} < \infty$ for some $\alpha>0$, then we can not apply Theorem \ref{theo:1}. However here $\delta_{n,0}\leq \big(\log (1+\sqrt{n})\big)^{-\alpha}$ and hence using theorem \ref{theo:3}, UCLT over $\A^{re}$ holds if $\log p = o\big((\log n)^{\alpha/3}\big)$. Therefore, $\log p$ can still grow sub-exponentially with $n$ even without existence of absolute polynomial moment  $\geq 2$. The growth in this case depends on how fast $\delta_{n, 0}$ goes to $0$ or in other words it depends on how fat the tail of the distribution of $X_{11}$ is.
 
 \begin{remark}
Our results of this section should be compared with the UCLT results available in the literature in case of the normalized sums. 
The rate $\log p =o\big(n^{1/5}\big)$ in the UCLT $\rho_{n,\A^{re}}\rightarrow 0$ of the normalized sum $U_n$ is established in \cite{CCKK} and \cite{K} when $\mathbf{E}\exp{\big(t|X_{ij}|\big)}< \infty$, for all $(i,j)$. On the other hand, here we show that to achieve $\log p =o(n^{1/5})$ in the UCLT $\tau_{n,\A^{re}}\rightarrow 0$ of the self-normalized sum $T_n$ we need $\mathbf{E}|X_{ij}|^{5/2}< \infty$, for all $(i,j)$.
Again the rate $\log p=o\big(n^{1/3}(\log n)^{-2/3}\big)$ in $\rho_{n,\A^{re}}\rightarrow 0$ is obtained by \cite{FK} under the assumption of the random vectors $X_i$'s having log-concave densities. Whereas for the rate $\log p=o\big(n^{1/3}\big)$ in $\tau_{n,\A^{re}}\rightarrow 0$, we need only the existence of third absolute moments of $X_{ij}$ for all $(i,j)$.  
 \end{remark}


\subsection{Rate of growth of $p$ under fourth moment} In this sub-section we try to improve the results obtained in the previous sub-section on the rate of growth of $p$ under higher moment conditions. In particular, we show that $\log p$ can grow like $o\big(n^{1/2}\big)$ under the existence of fourth moment. Additionally, the UCLT over $\A^{re}$ holds with an error $\epsilon$ whenever $\log p =\epsilon n^{1/2}$ for a non trivial set of $\epsilon>0$. Now the natural question is if the growth rate of $\log p$ can be improved further from $n^{1/2}$ as a power of $n$ under stronger moment conditions. The answer is \textit{no} even when $X_{ij}$'s are iid across $(i,j)$ and $X_{11}$ is stochastically bounded. The rate of growth of $\log p$ in the UCLT of $T_n$ over $\A^{re}$ essentially stabilizes at $n^{1/2}$ whenever fourth or higher absolute moments of $X_{ij}$ exists. We are now ready to state our first theorem of this sub-section. 

\begin{theo}\label{theo:4}
Let $\max_{i=1,\dots,n}\mathbf{E}|X_{i1}|^4 = O(1)$, $\liminf_{n\rightarrow \infty}\sigma_n^2 >0$ and $EX_{i1}^3 = 0$ for all $i \in \{1,\dots, n\}$. Then we have $$\tau_{n,\mathcal{A}^{re}}\rightarrow 0 \;\;\; \text{as}\;\;n\rightarrow \infty,$$
provided $\log p =o(n^{1/2})$. Moreover, there exists a positive constant $c\leq 1/8$ such that whenever $\log p \leq \epsilon n^{1/2}$ with $0<\epsilon < c$, $$\limsup_{n \rightarrow \infty}\tau_{n,\mathcal{A}^{re}}< \epsilon.$$
\end{theo}

The next proposition shows that for a UCLT
even over the smaller class of sets $\A^{max}$, the $o\big(n^{1/2}\big)$ upper 
bound on $\log p$ can not be significantly improved upon even for the smaller class $\mathcal{A}^{max}$.

\begin{proposition}\label{prop:2.2}
Let $X_{ij}$'s be iid Rademacher variables, i.e. $X_{ij}=1$ or $-1$ each with probability $1/2$ and be independent across $i\in \{1\dots,n\}$ and $j\in \{1,\dots, p\}$. If 
$\limsup_{\nti} \big[(n\log n)^{-1/2}\log p\big] >0$, then 
\begin{align}\label{eqn:negative}
\tau_{n,\mathcal{A}^{max}}\nrightarrow 0\;\; \text{as}\;\; n \rightarrow \infty.
\end{align}
\end{proposition}


\section{An application to Student's t statistic}\label{sec:student's t} In this section we apply the results obtained in the previous section for finding UCLT for high dimensional component-wise Student's t statistic. Recall that based on the mean zero random vectors $X_1,\dots,X_n$, the corresponding student's t statistic is $W_{n}=(W_{n1},\dots,W_{n,p})^\prime$, where $$W_{nj}=\frac{\sqrt{n}\bar{X}_{nj}}{\sqrt{(n-1)^{-1}\sum_{i=1}^{n}(X_{ij}-\bar{X}_{nj})^2}},$$
where $\bar{X}_{nj}=n^{-1}\sum_{i=1}^{n}X_{ij}$. The statistic $W_n$ can be used for drawing high dimensional inference based on the sample $X_1,\dots,X_n$. In fact in most of the practical applications, component-wise studentization is more natural to perform than component-wise standardization since component-wise population variances are not in general available. In statistical terms, $W_n$ is always a statistic irrespective of whether underlying variance structure is known or unknown. This is not true for the standardized sum $U_n$. 
To that end, we are interested in high dimensional UCLT for $W_n$ by investigating the quantity given by
\begin{align}
    \gamma_{n,\A}\equiv   \sup_{A\in \A} \Big| P(W_n\in A) - P(Z\in A)   \Big| 
    \label{eqn:clt-3}
\end{align}
for some class $\A$ of convex subsets of $\mathcal{R}^p$.    
We show that results similar to Section \ref{sec:main} continue to hold for $\gamma_{n,\A}$ whenever $\A=\A^{re}$ or $\A^{max}$ or $\A^{dist}$. Recall the relation (\ref{eqn:relation}) between $T_{nj}$ and $W_{nj}$ that is crucial in establishing a high dimensional UCLT for $W_n$ from that for $T_n$. 
We are now ready to state the high dimensional UCLT result for $W_n$.
\begin{theo}\label{theo:6}
Statements of the Theorems \ref{theo:1}, \ref{theo:3} and \ref{theo:4} are true if $\tau_{n,\mathcal{A}^{re}}$ is replaced by $\gamma_{n,\mathcal{A}^{re}}$.
\end{theo}
Theorem \ref{theo:6} shows that all the high dimensional UCLT results that we developed for self-normalized sums, also hold for component-wise Student's t-statistic. Therefore the same rate of $\log p$ can be achieved in the high dimensional UCLT over $\A^{re}$ under much weaker moment conditions if we use studentzation in place of standardization. \cite{FK} obtained the rate $\log p=o\big(n^{1/3}(\log n)^{-2/3}\big)$ in $\rho_{n,\A^{re}}\rightarrow 0$ under the assumption of the random vectors $X_i$'s having log-concave densities (i.e. when $X_{ij}$ are sub-exponential for all $(i,j)$). Whereas for the rate $\log p=o\big(n^{1/3}\big)$ in $\tau_{n,\A^{re}}\rightarrow 0$, we only need the existence of third absolute moments of $X_{ij}$ for all $(i,j)$. Again to achieve $\log p =o\big(n^{1/2}\big)$ in $\gamma_{n,\A^{re}}\rightarrow 0$ we need $\mathbf{E}|X_{ij}|^4< \infty$ and $\mathbf{E}X_{ij}^3 = 0$, for all $(i,j)$. Whereas $X_{ij}$'s are assumed to be sub-Gaussian to achieve $\log p =o\big(n^{1/2}\big)$ in $\rho_{n,\A^{re}}\rightarrow 0$, as established in \cite{DL}. 
Additionally under the more specialized iid structure of $X_{ij}$'s, sub-exponential growth rate of $p$ can be achieved in (\ref{eqn:clt-3}) with $\mathcal{A}=\mathcal{A}^{re}$ without any moment conditions, only the assumption of $X_{11}$ being in the domain of attraction of normal distribution is sufficient.

\section{Near classical Berry-Esseen Rate}\label{sec:BE}
In this section we establish near classical Berry-Esseen rate for both self-nomalized and studentized random vectors in high dimension. Let us consider following two cases separately for describing the results.\\
Case I: $X_{ij}$'s are all iid and $X_{11}$ is in the domain of attraction of the normal distribution.\\
Case II: $X_i =(X_{i1},\dots, X_{ip})^\prime$ are independent random vectors for $i\in \{1,\dots, n\}$. For each $i$, $X_{ij}$'s are iid, $EX_{i1}=0$ and $E|X_{i1}|^{2+\kappa} < \infty$ for some $0< \kappa \leq 1$.\\

Let us consider Case I first. Under this setup, \cite{BG} proved that $$\tau_{n,\mathcal{A}^{dist}} \leq A_1\delta_{n,0}$$ for some absolute positive constant $A_1$, when $p=1$ (cf. Theorem 1.4 of \cite{BG} and Theorem 3.1 of \cite{WS}), where $\delta_{n, 0}$ is defined in the beginning of Section \ref{sec:main}. Here we establish that $$\tau_{n,\mathcal{A}^{dist}} \leq M_1\delta_{n,0}\Big(\log\Big(\dfrac{p}{\delta_{n, 0}}\Big)\Big)^3,$$ for some positive constant $M_1$, even when $p$ grows with $n$. Clearly for fixed $p$, our bound nearly recovers the classical Berry-Esseen rate of \cite{BG}. The next theorem states the high dimensional Berry-Esseen rate obtained under the setup of Case I in more precise manner.

\begin{theo}\label{theo:7}
Let $X_{ij}$'s are iid copies of $X_{11}$ for all $i\in\{1,\dots,n\}$ and $j\in \{1,\dots,p\}$ where $p\geq 3$. Let $X_{11}$ be in the domain of attraction of the normal distribution with the distribution of $X_{11}$ being non-degenerate. 
Then there exists some positive constant $M_1$ such that
\begin{align*}
\tau_{n,\mathcal{A}^{re}}, \gamma_{n, \mathcal{A}^{re}} \leq M_1\delta_{n,0}\Big(\log\Big(\dfrac{p}{\delta_{n, 0}}\Big)\Big)^3
\end{align*}
\end{theo}
The constant $M_1$ in the above theorem can be made precise. Indeed it can be shown that $M_1$ can be taken as $1664C_3$ where the constant $C_3$ is same as that appearing in Lemma \ref{lem:3} but with $Z_1 = X_{11}$. Now let us consider Case II. Under this setup, \cite{S05} established that 
$$\tau_{n, dist}\leq \dfrac{25}{d_{n,\kappa}^{2+\kappa}},$$ when $\kappa \in (0,1]$ and $p=1$, where $d_{n, \kappa}$ is defined in the beginning of Section \ref{sec:main}. See also Theorem 1.2 of \cite{BBG} and Theorem 1 of \cite{N}. 
Here we establish a high dimensional Berry-Esseen bound with dependence on $p$ in terms of $(\log p )^{(2+\kappa)/2}$. 
Next theorem summarizes this finding.

\begin{theo}\label{theo:8}
Let $X_i =(X_{i1},\dots, X_{ip})^\prime$ be independent random vectors in $\mathcal{R}^p$ for $i\in \{1,\dots, n\}$ where for each $i$, $X_{ij}$'s are iid, $EX_{i1}=0$ and $E|X_{i1}|^{2+\kappa} < \infty$ for some $0< \kappa \leq 1$. Assume that $p\geq 3$ and $d_{n, \kappa} \geq 1$. Then there exists an absolute constant $A_{\kappa}$ (which depends on $\kappa$ only) such that
\begin{align*}
\tau_{n,\mathcal{A}^{re}}, \gamma_{n, \mathcal{A}^{re}}  \leq\; & \dfrac{A_{\kappa}(\log p d_{n,\kappa})^{(2+\kappa)/2}}{d_{n,\kappa}^{2+\kappa}}.
\end{align*}
\end{theo}

More precisely, the proof dictates that $A_{\kappa}$ can be taken as $A2^{10+\kappa}[2(2+\kappa)]^{(2+\kappa)/2}$ where the absolute constant $A$ is same as that is appearing in Lemma \ref{lem:3}. Now let us consider the case when $\kappa =1$ in Theorem \ref{theo:8} and assume that $n^{-1}\sum_{i=1}^{n}\mathbf{E}|X_{i1}|^{2+\kappa} = O(1)$ and $\liminf_{n\rightarrow \infty}\Big(n^{-1}\sum_{i=1}^{n}$ $EX_{i1}^2\Big) >0$. Then the rate obtained in Theorem \ref{theo:8} becomes
\begin{align*}
\dfrac{(\log p n)^{(2+\kappa)/2}}{\sqrt{n}}.
\end{align*}
Therefore, near-$n^{-1/2}$ convergence rate can be achieved in high dimensional UCLT for self-normalized/ studentized sums under the existence of third absolute moments. When $\kappa =1$, the above rate  can be compared with the rate obtained for high dimensional normalized sums in \cite{FK} which was obtained under the assumption of $X_i$'s being iid having log-concave density. 

The rate of growth of the dimension $p$ obtained in Theorem \ref{theo:8} is generally the optimal convergence rate for both self-nomalized and studentized sums under the existence of third absolute moments, i.e. when $\kappa = 1$ in Theorem \ref{theo:8}. This is stated in the next proposition for self-normalized sums for the simplest case of $X_{ij}$'s being iid across $(i, j)$. 

\begin{proposition}\label{prop:4.1}
Let $X_i =(X_{i1},\dots, X_{ip})^\prime$, $i \in \{1,\dots, n\}$, be random vectors in $\mathcal{R}^p$ where $X_{ij}$'s are all iid. Also assume that $EX_{11}=0$, $EX_{11}^2 =1$, $E|X_{11}|^3 < \infty$ and $EX_{11}^3 > 0$. Then if $\log p =o\big(n^{1/3}\big)$ and $\sqrt{n} = o\Big(p(\log p)^{3/2}\Big)$ as $n\rightarrow \infty$, then  
\begin{align*}
\liminf_{n\rightarrow \infty}\sqrt{\dfrac{n}{(\log p)^3}}\tau_{n, \mathcal{A}^{\max}}>0.
\end{align*}
\end{proposition}
Above proposition can be compared with Proposition 1.1 of \cite{FK} where similar result was obtained for normalized sums. 


\begin{remark}
Our proofs of the theorems mainly depend on regrouping the endpoints of the $p$ intervals of the elements of $\mathcal{A}^{re}$ into a specific partition of the real line (depending on n) and then applying a set of suitable error
estimates over each range using non-uniform Berry-Esseen bounds for one dimensional self-normalized sum (cf. Lemmas \ref{lem:3}, \ref{lem:5} and \ref{lem:4}). As pointed out by a referee, an alternative to this approach is to use coupling inequalities combined with sub-Gaussian tail bounds of one dimensional self normalized sums and the fact that coupling inequality for maximum lead to UCLTs (cf. Lemma 2.1 of \cite{CCK16}). Either the non-uniform Berry-Esseen bound for one dimensional self-normalized sums required in our proofs or the coupling inequalities (like the results in \cite{MZ}) required in the alternative strategy both are byproduct of the large deviation results of one dimensional self normalized sums (cf. \cite{JSW}). However, the sub-Gaussian tail bound of self-normalized sums (required in the alternative strategy) generally requires the assumption of $X_1,\dots, X_n$ being iid ( cf. Theorem 2.5 of \cite{GGM}) or the assumption of the components of $X_i$'s being symmetric about $0$ when $X_1,\dots, X_n$ are assumed to be independent only (cf. Theorem 2.1 of \cite{WJ}). Neither of these assumptions is required in our proof strategy and the theorems are established even when $X_i$'s are independent but not necessarily identically distributed. The alternative proof technique for Theorem \ref{theo:8} when $\kappa =1$ is presented in the supplementary material file, under the assumption that $X_{ij}$'s are iid across $(i, j)$.
\end{remark}

\section{Proofs of the Results}
\label{sec:pf}
\setcounter{equation}{0} 
Suppose, $\Phi(\cdot)$ and $\phi(\cdot)$ respectively denote the cdf and pdf of the standard normal random variable in any dimension. Define $N_i=(N_{i1},\dots,N_{ip})^\prime$, $i\in \{1,\dots,n\}$ where $N_{ij}$'s are iid $N(0,1)$ random variables for all $i\in \{1,\dots,n\}$ and $j\in \{1,\dots,p\}$. For any vector $\bm{t}
=(t_1,\dots,t_p)\in \mathcal{R}^p$, let $t_{(j)}$ and $t^{(j)}$ respectively denote the $j$th element after sorting the components of $\bm{t}$ 
 in increasing order and in decreasing order.
 (We use boldface font only for $\bm{t}$ to avoid some 
 notational conflict later on. All other vectors are 
 denoted using regular font). 
For any random variable $H$, $P\big(H\leq x\big)$ is assumed to be $1$ if $x=\infty$. Any absolute constant is denoted by $A$. $C,C_1,C_2,\dots$ denote generic constants which depend only on the underlying distribution of the random vectors $X_1,\dots,X_n$. 

Recall that for a sequence of random variables $\{Y_i\}_{i\geq 1}$, $\sigma_n^2(Y)=n^{-1}\sum_{i=1}^{n}Var(Y_i)$, for any $k> 0$, $\beta_{n,k}(Y)=n^{-1}\sum_{i=1}^{n}\mathbf{E}|Y_i|^{2+k}$ and the Lyapunov's ratio is $$d_{n,k}(Y)=\dfrac{\sqrt{n}\sigma_n(Y)}{\big(n\beta_{n,k}(Y)\big)^{1/{(2+k)}}}.$$   
For some random variable $Y$ recall that
\begin{align*}
\delta_{n,x}(Y) = \;&n\mathbf{P}\Big(|Y|>\kappa_{n,x}(Y)\Big)+n\big[\kappa_{n,x}(Y)\big]^{-1}\Big|\mathbf{E}\big\{YI\big(|Y|\leq \kappa_{n,x}(Y)\big)\big\}\Big|\\
&+n\big[\kappa_{n,x}(Y)\big]^{-3}\mathbf{E}\big\{|Y|^3I\big(|Y|\leq \kappa_{n,x}(Y)\big)\big\},
\end{align*}
where $\kappa_{n,x}(Y)=\sup\Big\{s:ns^{-2}\mathbf{E}\big\{Y^2I(|Y|\leq s)\big\}\geq 1+x^2\Big\}$. Whenever $X_{ij}$'s are the underlying random variables, then we simply write $\sigma_n^2$, $\beta_{n,k}$, $d_{n,k}$ and $\kappa_{n,x}$. 

We will need some lemmas which are stated next. Proofs of the theorems
are given in Section 4.2 below.

\subsection{Auxiliary Lemmas}
\begin{lem}\label{lem:1}
For any $t>0$, $\dfrac{1}{t}\geq\dfrac{1-\Phi(t)}{\phi(t)}\geq \dfrac{2}{\sqrt{t^2+4}+t}$.
\end{lem}

\noindent
{\bf Proof of Lemma \ref{lem:1}:}
This inequality is proved in Birnbaum (1942).\\



\begin{lem}\label{lem:3}
Let $Z_1\dots,Z_n$ be independent random variables with $E|Z_i|^{2+\delta}<\infty$ for all $i\in\{1,\dots,n\}$ for some $0<\delta\leq 1$. Then for $0\leq |x| < d_{n,\delta}(Z)$, we have 
\begin{align*}
\Big|\mathbf{P}\Big(\frac{\sum_{i=1}^{n}Z_i}{\sqrt{\sum_{i=1}^{n}Z_i^2}}\leq x\Big)-\Phi(x)\Big|\leq A (1+|x|)^{1+\delta}e^{-x^2/2}\big[d_{n,\delta}(Z)\big]^{-(2+\delta)},    
\end{align*}
where $A>0$ is an absolute constant.
\end{lem}

\noindent
{\bf Proof of Lemma \ref{lem:3}:} 
This lemma follows from (2.8) and (2.9) of Theorem 2.3 of Jing et al. (2003) and applying Lemma \ref{lem:1}.

\begin{lem}\label{lem:5}
Suppose that $Z,Z_1,\dots,Z_n$ are iid non-degenerate random variables. $Z$ is in the domain of attraction of the normal law. Recall the definition of $\delta_{n,0}(Z)$ from Section \ref{sec:main}. Define, $\omega_n(Z) = \big[\delta_{n,0}(Z)\big]^{-1/6}$. Then there exists a constant $C_3>0$, independent of $n,x$, such that for $n\geq C_3$, 
\begin{align}\label{eqn:donc0}
\Big|\mathbf{P}\Big(\dfrac{\sum_{i=1}^{n}Z_i}{\sqrt{\sum_{i=1}^{n}Z_i^2}}\leq x\Big)-\Phi(x)\Big|\leq C_3 \delta_{n,0}(Z)\max\{1,x^5\}e^{-x^2/2},
\end{align}
whenever $ 2|x|\leq \omega_n(Z)$.
\end{lem}
{\bf Proof of Lemma \ref{lem:5}:} When $2|x|\leq 2$, then (\ref{eqn:donc0}) follows directly from Theorem 1.4 of Bentkus and G\"otze (1996). Now let us look into when $2 < 2|x|\leq \omega_n(Z)$. Recall the definitions of $\kappa_{n,x}(Z)$ and $\delta_{n,x}(Z)$ from Section \ref{sec:main}. Then from Remark 3 of Robinson and Wang (2005), we have
\begin{align}\label{eqn:donc1}
 \Big|\mathbf{P}\Big(\dfrac{\sum_{i=1}^{n}Z_i}{\sqrt{\sum_{i=1}^{n}Z_i^2}}\leq x\Big)-\Phi(x)\Big|\leq A\delta_{n,|x|}(Z)(1+|x|)^{-1}e^{-x^2/2}, 
\end{align}
for some absolute constant $A>0$, whenever $\delta_{n,|x|}(Z)\leq 1$. Now the aim is to replace $\delta_{n, |x|}(Z)(1+|x|)^{-1}$ in RHS of (\ref{eqn:donc1}) with $\delta_{n,0}(Z)|x|^5$ when $2\leq 2|x|\leq \omega_n(Z)$. To do that we need both upper and lower bound on $\kappa_{n,|x|}(Z)$ in terms of $\kappa_{n,0}(Z)$. It is easy to see that $\kappa_{n,|x|}(Z)\leq \kappa_{n,0}(Z)$ for any $x$. Regarding the lower bound, we claim that  $\kappa_{n,0}(Z)/(1+x^2)\leq \kappa_{n,|x|}(Z)$ for large enough $n$ uniformly for $2\leq 2|x|\leq \omega_n(Z)$.
Note that using Lemma 1.3 of Bentkus and G\"otze (1996) we have \begin{align}\label{eqn:donc2}
\frac{\kappa_{n,0}(Z)}{1+x^2}=\sup\Big\{s:ns^{-2}\mathbf{E}Z^2I\big(|Z|\leq s(1+x^2)\big) = 1+x^2\Big\} 
\end{align}
and $k_{n,0}(Z)/\sqrt{n}$ to be positive for sufficiently large $n$. Hence, $k_{n,0}(Z)/(1+x^2) \rightarrow \infty$ as $n\rightarrow \infty$, uniformly for $|x|\leq \omega_n(Z)$ since $\omega_n(Z)=o(n^{1/4})$. To show $\kappa_{n,0}(Z)/(1+x^2)\leq \kappa_{n,|x|}(Z)$  uniformly in $2\leq 2|x|\leq \omega_n(Z)$ for sufficiently large $n$, enough to show $\mathbf{E}Z^2I\big(|Z|\leq k_{n,0}(Z)\big)\leq (1+x^2)\mathbf{E}Z^2I\Big(|Z|\leq k_{n,0}(Z)/(1+x^2)\Big)$, uniformly in $2\leq 2|x|\leq \omega_n(Z)$, for sufficiently large $n$. This is indeed true since $\mathbf{E}Z^2I\big(|Z|\leq s\big)$ is slowly varying as a function of $s$ at infinity [cf. Lemma 1(a) of C\"org\H{o} et al. (2003)] and the fact that $k_{n,0}(Z)/(1+x^2) \rightarrow \infty$ as $n\rightarrow \infty$ uniformly for $|x|\leq \omega_n(Z)$.    Therefore we have for sufficiently large $n$,
\begin{align}\label{eqn:donc3}
\delta_{n,|x|}(Z) = \;&n\mathbf{P}\Big(|Z|>\kappa_{n,|x|}(Z)\Big)+n\big[\kappa_{n,|x|}(Z)\big]^{-3}\mathbf{E}\big\{|Z|^3I\big(|Z|\leq \kappa_{n,|x|}(Z)\big)\big\}\nonumber\\
\leq\; & n\mathbf{P}\Big(|Z|>\kappa_{n,0}(Z)/(1+x^2)\Big)\nonumber\\
&+n(1+x^2)^3\big[\kappa_{n,0}(Z)\big]^{-3}\mathbf{E}\big\{|Z|^3I\big(|Z|\leq \kappa_{n,0}(Z)\big)\big\}\nonumber\\
\leq \; & 3(1+x^2)^3\delta_{n,0}(Z),
\end{align}
uniformly in $2\leq 2|x|\leq \omega_n(Z)$.
From (\ref{eqn:donc1}) and (\ref{eqn:donc3}), for sufficiently large $n$ we have
\begin{align}\label{eqn:donc4}
 \Big|\mathbf{P}\Big(\dfrac{\sum_{i=1}^{n}Z_i}{\sqrt{\sum_{i=1}^{n}Z_i^2}}\leq x\Big)-\Phi(x)\Big|&\leq 3A(1+x^2)^3\delta_{n,0}(Z)(1+|x|)^{-1}e^{-x^2/2}\nonumber\\
 & \leq 96A\delta_{n,0}(Z)|x|^5e^{-x^2/2},
\end{align}
whenever $x\in \Big\{\big[3(1+y^2)^3\delta_{n,0}(Z)\big]\leq 1\Big\}\cap \Big\{2\leq 2|y|\leq \omega_n(Z)\Big\}$. Again apply Lemma 1.3 of Bentkus and G\"otze (1996) and Lemma 1 of Cs\"org\H{o} et al. (2003) to claim that $\delta_{n,0}(Z)\rightarrow 0$ as $n\rightarrow \infty$. Therefore, (\ref{eqn:donc0}) follows from (\ref{eqn:donc4}). 

\begin{lem}\label{lem:4}
Let $Z_1\dots,Z_n$ be independent random variables such that $\max_{i=1,\dots,n}$ $E|Z_i|^{4}=O(1)$, $\liminf_{n\rightarrow \infty} \sigma_n^2(Z) >0$ and $\sum_{i=1}^{n}EZ_i^3 =O (n^{3/4})$. Then whenever $|x|\leq n^{1/4}$, for sufficiently large $n$ we have
\begin{align*}
\Big|\mathbf{P}\Big(\frac{\sum_{i=1}^{n}Z_i}{\sqrt{\sum_{i=1}^{n}Z_i^2}}\leq x\Big)-\Phi(x)\Big|\leq C_2 \big(1+|x|\big)e^{-x^2/2}n^{-1/2},    
\end{align*}
for some constant $C_2>0$, independent of $n,x$.
\end{lem}

\noindent
{\bf Proof of Lemma \ref{lem:4}:} 
When $0\leq x\leq n^{1/4}$, the statement is a direct consequence of Corollary 2.2 of Sang and Ge (2017). The statement for $0\geq x\geq -n^{1/4}$ follows by replacing $Z_i$ by $-Z_i$ and then applying Corollary 2.2 of Sang and Ge (2017).

\begin{lem}\label{lem:6}
For any positive integer $m$, $$\sqrt{2\pi}\;m^{m+1/2}e^{-m}\leq m! \leq  m^{m+1/2}e^{-m+1}.$$
\end{lem}
This is the well-known Stirling's formula. See for example Robbins (1955).

\begin{lem}\label{lem:7}
Let $Z_1,\dots,Z_n$ be a sequence of mean zero independent random vectors in $\mathcal{R}^p$ with $Z_i = (Z_{i1},\dots Z_{ip})$, $i \in \{1,\dots, n\}$ and let  $\{Z_{i1},\dots,Z_{ip}\}$ be
iid for each $i \in \{1,\dots, n\}$ with $\tilde{d}_{n,j}^2= \sum_{i=1}^{n}Z_{ij}^2$. Define, $l_1(x)=\max\Big{\{}\mathbf{P}\Big(\tilde{d}_{n1}^{-1}\sum_{i=1}^{n}\big(-Z_{i1}\big) \leq x \Big), \Phi(x)\Big{\}}$,  $d_1(x)=\Big{|}\mathbf{P}\Big(\tilde{d}_{n,1}^{-1}\sum_{i=1}^{n}\big(-Z_{i1}\big) \leq x \Big)-\Phi(x)\Big{|}$, $l_2(x)=\max\Big{\{}\mathbf{P}\Big(\tilde{d}_{n1}^{-1}$ $\sum_{i=1}^{n}Z_{i1} \leq x \Big), \Phi(x)\Big{\}}$ and  $d_2(x)=\Big{|}\mathbf{P}\Big(\tilde{d}_{n,1}^{-1}\sum_{i=1}^{n}Z_{i1} \leq x \Big)-\Phi(x)\Big{|}$. Also define, $\tilde{D}_n=diag\big(\tilde{d}_{n,1},\dots,\tilde{d}_{n,p}\big)$. Then we have
\begin{align*}
&\Big{|}\mathbf{P}\Big(\tilde{D}_n^{-1}\sum_{i=1}^{n}Z_{i} \in \prod_{j=1}^{p}\big{\{}[a_j,b_j]\cap \mathcal{R}\big{\}}\Big)-\mathbf{P}\Big(n^{-1/2}\sum_{i=1}^{n}N_{i} \in \prod_{j=1}^{p}\big{\{}[a_j,b_j]\cap\mathcal{R}\big{\}}\Big)\Big{|}\\
&\;\leq L_1(\bm{a}) +L_2(\bm{b}),
\end{align*}
where $\bm{a}=(a_1,\dots,a_p)^\prime$, $\bm{b}=(b_1,\dots,b_p)^\prime$, $$L_1(\bm{a})=\bigg[\sum_{k=1}^{p}\Big(\prod_{j\neq k}l_1\big(-a^{(j)}\big)\Big)d_1\big(-a^{(k)}\big)\bigg],\;\;\;\; L_2(\bm{b})=\bigg[\sum_{k=1}^{p}\Big(\prod_{j\neq k}l_2\big(b_{(j)}\big)\Big)d_2\big(b_{(k)}\big)\bigg].$$ 
\end{lem}

\noindent
{\bf Proof of Lemma \ref{lem:7}:}
Note that $\tilde{D}_n^{-1}\sum_{i=1}^{n}Z_{i}=(W_1,\ldots, W_p)'$
where $W_j=\tilde{d}_{n,j}^{-1}\sum_{i=1}^{n}$ $Z_{ij}$, $j\in \{1,\dots,p\}$, Then, using the nature of $Z_{ij}$'s, it is easy to see that $\{W_1,\dots,W_p\}$ are iid. Similarly, since $N_{ij}\sim N(0,1)$ are iid, the 
$p$-variables $\Big(n^{-1/2}\sum_{i=1}^{n}N_{i1}\Big),\dots,$ $\Big(n^{-1/2}\sum_{i=1}^{n}N_{ip}\Big)$ are also  iid. Hence we have
\begin{align*}
&\Big{|}\mathbf{P}\Big(\tilde{D}_n^{-1}\sum_{i=1}^{n}Z_{i} \in \prod_{j=1}^{p}\big{\{}[a_j,b_j]\cap \mathcal{R}\big{\}}\Big)-\mathbf{P}\Big(n^{-1/2}\sum_{i=1}^{n}N_{i} \in \prod_{j=1}^{p}\big{\{}[a_j,b_j]\cap\mathcal{R}\big{\}}\Big)\Big{|}\\
& = \bigg|\prod_{j=1}^{p}P\Big(W_j \in [a_j,b_j]\cap \mathcal{R}\Big)-\prod_{j=1}^{p}P\Big(n^{-1/2}\sum_{i=1}^{n}N_{ij} \in [a_j,b_j]\cap \mathcal{R}\Big)\bigg|\\
& \leq \bigg[\sum_{k=1}^{p}\Big(\prod_{j\neq k}\Big(\min\Big{\{}l_1\big(-a_j\big),l_2\big(b_j\big)\Big{\}}\Big)\Big)\Big[d_1\big(-a_k\big)+d_2\big(b_k\big)\Big]\bigg]\\
& \leq \sum_{k=1}^{p}\Big(\prod_{j\neq k}l_1\big(-a_j\big)\Big)\Big[d_1\big(-a_k\big)\Big]+ \sum_{k=1}^{p}\Big(\prod_{j\neq k}l_2\big(b_j\big)\Big)\Big[d_2\big(b_k\big)\Big]\\
& = \sum_{k=1}^{p}\Big(\prod_{j\neq k}l_1\big(-a^{(j)}\big)\Big)\Big[d_1\big(-a^{(k)}\big)\Big]+ \sum_{k=1}^{p}\Big(\prod_{j\neq k}l_2\big(b_{(j)}\big)\Big)\Big[d_2\big(b_{(k)}\big)\Big]
\end{align*}
The last equality is due to the following fact:\\
If $(G_1,H_1),\dots,(G_p,H_p)$ are iid random vectors in $\mathcal{R}^2$, then for any $t_1\dots,t_p \in \mathcal{R}$,
\begin{align*}
&\sum_{k=1}^{p}\bigg[\Big(\prod_{j \neq k}\Big(\max\Big{\{}P\Big(G_j \leq t_j\Big), P\Big(H_j\leq t_j\Big)\Big{\}}\Big)\Big)\Big|P\Big(G_k\leq t_k\Big)-P\Big(H_k\leq t_k\Big)\Big|\bigg]\\
&=
\sum_{k=1}^{p}\bigg[\Big(\prod_{j \neq k}\Big(\max\Big{\{}P\Big(G_1 \leq t_{(j)}\Big), P\Big(H_1\leq t_{(j)}\Big)\Big{\}}\Big)\Big)\Big|P\Big(G_1\leq t_{(k)}\Big)-P\Big(H_1\leq t_{(k)}\Big)\Big|\bigg],
\end{align*}
where $\{t_{(1)}, t_{(2)},\dots, t_{(p)}\}$ are obtained after sorting $\{t_1,\dots,t_p\}$ in increasing order.
Therefore we are done.

\subsection{Proofs of the main results}

\noindent
{\bf Proof of Theorem \ref{theo:1}:}  We are going to prove part (b) only. Part (a) follows from part (b) by taking $\epsilon\rightarrow 0$. Recall that $T_n=\Big(\frac{\sum_{i=1}^{n}X_{i1}}{\sqrt{\sum_{i=1}^{n}X_{i1}^2}},\dots,\frac{\sum_{i=1}^{n}X_{ip}}{\sqrt{\sum_{i=1}^{n}X_{ip}^2}}\Big)^\prime$ and suppose that $S_n=n^{-1/2}\sum_{i=1}^{n}N_i$. Let $T=(T_{n1},\dots,T_{np})^\prime$ and $S_n= (S_{n1},\dots,S_{np})^\prime$. Clearly $T_{nj}$'s are iid and $S_{nj}$'s are iid for $j\in \{1,\dots,p\}$. We can use Lemma \ref{lem:7} with $Z_i=X_i$ for $i\in \{1,\dots,n\}$, to obtain
\begin{align*}
\Big{|}\mathbf{P}\Big(T_n \in \prod_{j=1}^{p}\big{\{}[a_j,b_j]\cap \mathcal{R}\big{\}}\Big)-\mathbf{P}\Big(S_n \in \prod_{j=1}^{p}\big{\{}[a_j,b_j]\cap\mathcal{R}\big{\}}\Big)\Big{|}
\leq  L_1(\bm{a}) +L_2(\bm{b}),
\end{align*}
where $L_1(\bm{a})$ and $L_2(\bm{b})$ are as defined in Lemma \ref{lem:7}. Since all the assumptions are also satisfied if we replace $\{X_1,\dots,X_n\}$ by $\{-X_1,\dots,-X_n\}$, it is enough to show
\begin{align}\label{eqn:4}
\sup_{t_1\leq t_2\leq \dots \leq t_p} L((t_1,\dots,t_p)^\prime)= \sup_{t_1\leq t_2\leq \dots \leq t_p}\bigg[\sum_{j=1}^{p}\Big(\prod_{j\neq k}l(t_j)\Big)d(t_k)\bigg]\leq \epsilon^{(1+\delta)/3}/2\;\; \text{for sufficiently large}\; n.
\end{align}
Here, $l(x)=\max\Big{\{}\mathbf{P}\Big(T_{n1} \leq x \Big), \mathbf{P}\Big(S_{n1} \leq x \Big)\Big{\}} \;\;\; \text{and}\;\;\; d(x)=\Big{|}\mathbf{P}\Big(T_{n1} \leq x \Big)-\mathbf{P}\Big(S_{n1} \leq x\Big)\Big{|}$. 

Now take $8c=\min\Big\{[A2^{1+\delta}\sqrt{2\pi}(\sqrt{2}+1)]^{-3},1\Big\}$. Fix $\bm{t}=(t_1,\dots,t_p)^\prime$ in $\mathcal{R}^p$ such that $t_1\leq t_2 \leq \dots\leq t_p$. Then there exist integers $l_1,l_2,l_3$, depending on $n$, such that $0\leq l_1, l_2,l_3 \leq p$ and
\begin{align}\label{eqn:14}
&t_1\leq t_2 \leq \dots \leq t_{l_1} < -\epsilon^{1/3} d_{n,\delta}\nonumber\\
-\epsilon^{1/3}d_{n,\delta} \leq\; & t_{l_1+1}\leq t_{l_1+2}\leq \dots \leq\; t_{l_2} < 1\nonumber\\
1 \leq\; & t_{l_2+1}\leq  t_{l_2+2}\leq \dots \leq t_{l_3} \leq \epsilon^{1/3} d_{n,\delta}\nonumber\\
\epsilon^{1/3}d_{n,\delta} <\; & t_{l_3+1}\leq  t_{l_3+2}\leq \dots \leq t_p 
\end{align}

Now use the same definitions of $l(x)$ and $d(x)$, as in the proof of Theorem \ref{theo:1}. Then due to Lemma \ref{lem:1}, \ref{lem:3} and the fact that $\epsilon \leq c < 1$, we have for sufficiently large $n$,
\begin{align}\label{eqn:15}
&l(x) \leq\; I\Big(x>\epsilon^{1/3}d_{n,\delta}\Big) + \bigg[1-\dfrac{2\phi(1)}{\sqrt{5}+1}+ A2^{1+\delta} e^{-1/2}d_{n,\delta}^{-(2+\delta)}\bigg]I\Big(x < 1\Big)\nonumber\\
&\;\;\;\;\;\;\;\;\;\;\;\; +\bigg[1-\dfrac{2\phi(x)}{\sqrt{x^2+4}+x}+A2^{1+\delta}\epsilon^{(1+\delta)/3}d_{n,\delta}^{-1} e^{-x^2/2}\bigg]I\Big(x \in\Big[1,\epsilon^{1/3}d_{n,\delta}\Big]\Big)\nonumber\\
\leq\; & I\Big(x>\epsilon^{1/3}d_{n,\delta}\Big) + \bigg[1-\dfrac{\phi(1)}{\sqrt{5}+1}\bigg]I\Big(x < 1\Big)+\bigg[1-\dfrac{d_{n,\delta}^{-1}e^{-x^2/2}}{\sqrt{2\pi}\Big(\sqrt{\epsilon^{2/3}+4d_{n,\delta}^{-2}}+\epsilon^{1/3}\Big)} \bigg]I\Big(x \in\Big[1,\epsilon^{1/3}d_{n,\delta}\Big]\Big),
\end{align}
and
\begin{align}\label{eqn:16}
d(x) \leq \Big[2e^{-\epsilon^{2/3}d_{n,\delta}^{2}}\Big]I\Big(|x|>\epsilon^{1/3}d_{n,\delta}\Big) + \Big[A2^{1+\delta}\epsilon^{(1+\delta)/3} d_{n,\delta}^{-1}e^{-x^2/2}\Big]I\Big(|x|\leq \epsilon^{1/3}d_{n,\delta}\Big)
\end{align}
for any $x\in \mathcal{R}$, for sufficiently large $n$. Therefore from equations (\ref{eqn:4})-(\ref{eqn:16}) we have for sufficiently large $n$,
\begin{align}\label{eqn:17}
L(\bm{t})\leq J_1(\bm{t}) +J_2(\bm{t}) +J_3(\bm{t}) + J_4(\bm{t}),
\end{align}
where 
\begin{align*}
J_1(\bm{t}) =\; & \bigg(\Big[1-\dfrac{\phi(1)}{\sqrt{5}+1}\Big]^{l_2-1}\bigg)*\bigg(\prod_{j=l_2+1}^{l_3}\bigg[1-\dfrac{(2\pi)^{-1/2}d_{n,\delta}^{-1}e^{-x^2/2}}{\Big(\sqrt{\epsilon^{2/3}+4d_{n,\delta}^{-2}}+\epsilon^{1/3}\Big)} \bigg]\bigg)*\bigg(\sum_{k=1}^{l_1}2e^{-\epsilon^{2/3}d_{n,\delta}^{2}/2}\bigg),
\end{align*}
\begin{align*}
J_2(\bm{t})=\;&\bigg(\Big[1-\dfrac{\phi(1)}{\sqrt{5}+1}\Big]^{l_2-1}\bigg)*\bigg(\prod_{j=l_2+1}^{l_3}\bigg[1-\dfrac{(2\pi)^{-1/2}d_{n,\delta}^{-1}e^{-x^2/2}}{\Big(\sqrt{\epsilon^{2/3}+4d_{n,\delta}^{-2}}+\epsilon^{1/3}\Big)} \bigg]\bigg)*\bigg(\sum_{k=l_1+1}^{l_2}\dfrac{A2^{1+\delta}\epsilon^{(1+\delta)/3}}{d_{n,\delta}e^{t_k^2/2}}\bigg),
\end{align*}
\begin{align*}
J_3(\bm{t})=\;&\bigg(\Big[1-\dfrac{\phi(1)}{\sqrt{2}+1}\Big]^{l_2}\bigg)*\bigg(\sum_{k=l_2+1}^{l_3}\dfrac{A2^{1+\delta}\epsilon^{(1+\delta)/3}}{d_{n,\delta}e^{t_k^2/2}}\bigg(\prod_{ \substack{j=l_2+1\\ j \neq k}}^{l_3}\bigg[1-\dfrac{d_{n,\delta}^{-1}e^{-x^2/2}}{\sqrt{2\pi}\Big(\sqrt{\epsilon^{2/3}+4d_{n,\delta}^{-2}}+\epsilon^{1/3}\Big)} \bigg]\bigg)\bigg),
\end{align*}
\begin{align*}
J_4(\bm{t})=\;&\bigg(\Big[1-\dfrac{\phi(1)}{\sqrt{5}+1}\Big]^{l_2}\bigg)*\bigg(\prod_{j=l_2+1}^{l_3}\bigg[1-\dfrac{d_{n,\delta}^{-1}e^{-x^2/2}}{\sqrt{2\pi}\Big(\sqrt{\epsilon^{2/3}+4d_{n,\delta}^{-2}}+\epsilon^{1/3}\Big)} \bigg]\bigg)*\bigg(\sum_{k=l_3+1}^{p}2e^{-\epsilon^{2/3}d_{n,\delta}^{2}/2}\bigg).
\end{align*}

\emph{Bound on $J_1(\bm{t})+J_4(\bm{t})$}:
Since $\log p =\epsilon d_{n,\delta}^2$, from (\ref{eqn:17}) we have
\begin{align}\label{eqn:18}
J_1(\bm{t})+J_4(\bm{t}) 
\leq p\Big(2e^{-\epsilon^{2/3}d_{n,\delta}^{2}}\Big)=2\exp{\big(\epsilon d_{n,\delta}^2-\epsilon^{2/3}d_{n,\delta}^2/2\big)}
< \epsilon^{(1+\delta)/3}/12,
\end{align}
for large enough $n$, since $\epsilon < c \leq 1/8$.

\emph{Bound on $J_2(\bm{t})$}:
Noting that $d^{-1}=\Big[1-\dfrac{\phi(1)}{\sqrt{5}+1}\Big]$ and $\epsilon<1$, we have for sufficiently large $n$,
\begin{align}\label{eqn:19}
J_2(\bm{t})\leq I_2(\bm{t}) &\leq \bigg(\Big[1-\dfrac{\phi(1)}{\sqrt{5}+1}\Big]^{l_2-1}\bigg)\bigg(\sum_{k=l_1+1}^{l_2}A2^{1+\delta}\epsilon^{(1+\delta)/3}d_{n,\delta}^{-1}e^{-t_k^2/2}\bigg)I\Big((l_2-l_1) \geq 1\Big)\nonumber\\
&\leq  \Big(d^{-1} (\log d)^{-1}d^{-(\log d)^{-1}}\Big) A2^{1+\delta}\epsilon^{(1+\delta)/3}d_{n,\delta}^{-1} < \epsilon^{(1+\delta)/3}/12.
\end{align}

\emph{Bound on $J_3(\bm{t})$}:
Note that if $(l_3-l_2) = 0$ then $J_3(\bm{t})=0$ and there is nothing more to do. Hence assume $(l_3-l_2)\geq 1$. Then we have
\begin{align*}
\dfrac{\partial J_{3}(\bm{t})}{\partial t_l }=& \bigg[A2^{1+\delta}\epsilon^{(1+\delta)/3}d_{n,\delta}^{-1}t_le^{-t_l^2/2}n^{-1/4}\prod_{\substack{j=l_2+1\\ j\neq l}}^{l_3}\bigg[1-\dfrac{d_{n,\delta}^{-1}e^{-t_j^2/2}}{\sqrt{2\pi}\Big(\sqrt{\epsilon^{2/3}+4d_{n,\delta}^{-2}}+\epsilon^{1/3}\Big)} \bigg]\bigg]\times\\
&\Bigg[\sum_{\substack{k=l_2+1\\k\neq l}}^{l_3}\Bigg(\bigg[1-\dfrac{d_{n,\delta}^{-1}e^{-t_k^2/2}}{\sqrt{2\pi}\Big(\sqrt{\epsilon^{2/3}+4d_{n,\delta}^{-2}}+\epsilon^{1/3}\Big)} \bigg]^{-1}\dfrac{d_{n,\delta}^{-1}e^{-t_k^2/2}}{\sqrt{2\pi}\Big(\sqrt{\epsilon^{2/3}+4d_{n,\delta}^{-2}}+\epsilon^{1/3}\Big)}\Bigg)-1\Bigg]
\end{align*}
Hence for any $l=l_2+1,\dots,l_3$, $\dfrac{\partial J_{3}(\bm{t})}{\partial t_l }\gtreqless 0$ if and only if 
\begin{align}\label{eqn:12}
\sum_{\substack{j=l_2+1 \\ j\neq l}}^{l_3}\dfrac{z_j}{1-z_j}\gtreqless 1,
\end{align}
where $z_j = d_{n,\delta}^{-1}e^{-t_j^2/2}\Big[\sqrt{2\pi}\Big(\sqrt{\epsilon^{2/3}+4d_{n,\delta}^{-2}}+\epsilon^{1/3}\Big)\Big]^{-1}$. Note that since $1 \leq t_{l_2+1}\leq \dots \leq t_{l_3}$, $1> z_{l_2+1} \geq \dots \geq z_{l_3} >0$ for sufficiently large $n$. Hence for sufficiently large $n$, $$\dfrac{z_{l_2+1}}{1-z_{l_2+1}} \geq \dots \geq \dfrac{z_{l_3}}{1-z_{l_3}},$$ due to the fact that $z/(1-z)$ is increasing for $z\in (0,1)$. Therefore from (\ref{eqn:12}) we can say that $I_{31}(\bm{t})$ is non-increasing in $\{t_{l_2+1},\dots,t_m\}$ and non-decreasing in $\{t_{m+1},\dots,t_{l_3}\}$ where $(m-l_2)$ is a non-negative integer not more than $(l_3-l_2)$. Clearly $m$ is a function of $\bm{t}$. Write $\tilde{z}_n = d_{n,\delta}^{-1}e^{-1/2}\Big[\sqrt{2\pi}\Big(\sqrt{\epsilon^{2/3}+4d_{n,\delta}^{-2}}+\epsilon^{1/3}\Big)\Big]^{-1}$. Then for sufficiently large $n$, we have
\begin{align}\label{eqn:20}
J_3(\bm{t}) \leq\;& (q-l_2)\Big[1-\tilde{z}_n \Big]^{m-l_2-1}\Big(A2^{1+\delta}\epsilon^{(1+\delta)/3}d_{n,\delta}^{-1}e^{-1/2}\Big) \nonumber\\
&+ (l_3-m)\Big(A2^{1+\delta}\epsilon^{(1+\delta)/3}d_{n,\delta}^{-1}e^{-(2^{-1}\epsilon^{2/3}d_{n,\delta}^2)}\Big)\nonumber\\
\leq\;& 2\Big[\exp\Big(\log (m-l_2)-(m-l_2)\tilde{z}_n\Big)\Big]\Big(A2^{1+\delta}\epsilon^{(1+\delta)/3}d_{n,\delta}^{-1}e^{-1/2}\Big)\nonumber\\
&+ p\Big(A2^{1+\delta}\epsilon^{(1+\delta)/3}d_{n,\delta}^{-1}e^{-(2^{-1}\epsilon^{2/3}d_{n,\delta}^2)}\Big)\nonumber\\
\leq\;&  2\Big[\exp\Big(\sup_{x>0}\big[\log x-x\tilde{z}_n\big]\Big)\Big]\Big(A2^{1+\delta}\epsilon^{(1+\delta)/3}d_{n,\delta}^{-1}e^{-1/2}\Big)\nonumber\\
&+ \Big(A2^{1+\delta}\epsilon^{(1+\delta)/3}d_{n,\delta}^{-1}\exp{\big(-\epsilon^{2/3}d_{n,\delta}^2\big(1/2-\epsilon^{1/3}\big)}\Big)\nonumber\\
\leq\;& 2\tilde{z}_n^{-1}\Big(A2^{1+\delta}\epsilon^{(1+\delta)/3}d_{n,\delta}^{-1}e^{-1/2}\Big) + \Big(A2^{1+\delta}\epsilon^{(1+\delta)/3}d_{n,\delta}^{-1}\exp{\big(-\epsilon^{2/3}d_{n,\delta}^2\big(1/2-\epsilon^{1/3}\big)}\Big)\nonumber\\
\leq\;& \Big[\sqrt{2\pi}\big(\sqrt{\epsilon^{2/3}+4d_{n,\delta}^{-2}}+\epsilon^{1/3}\big)\Big]A2^{2+\delta}\epsilon^{(1+\delta)/3}\nonumber\\
&+ \Big(A2^{1+\delta}\epsilon^{(1+\delta)/3}d_{n,\delta}^{-1}\exp{\big(-\epsilon^{2/3}d_{n,\delta}^2\big(1/2-\epsilon^{1/3}\big)}\Big)\nonumber\\
 <\;& \epsilon^{(1+\delta)/3}/4 +\epsilon^{(1+\delta)/3}/12,
\end{align}
since $\epsilon^{1/3} < c^{1/3}\leq [A2^{4+\delta}\sqrt{2\pi}(\sqrt{2}+1)]^{-1}$. Now combining (\ref{eqn:14})-(\ref{eqn:20}), the proof of the part (b) of Theorem \ref{theo:1} is complete.\\

\noindent
{\bf Proof of Theorem \ref{theo:3}:}  We are going to prove part (b) only. Part (a) follows from part (b) by taking $\epsilon\rightarrow 0$. We are going to follow the same steps as in the proof of Theorem \ref{theo:1} with different estimates of $l(x)$ and $d(x)$. 

Now take $c=\min\Big\{\big[4C_3\sqrt{2\pi}(\sqrt{2}+1)\big]^{-3},1/8\Big\}$ where $C_3$ is the constant defined in Lemma \ref{lem:5} with $Z$ replaced by $X_{11}$. Recall that $\omega_n=\delta_n^{-1/6}$ where $\delta_n=\delta_{n,0}\big(X_{11}\big)$. Fix $\bm{t}=(t_1,\dots,t_p)^\prime$ in $\mathcal{R}^p$ such that $t_1\leq t_2 \leq \dots\leq t_p$. Then there exist integers $l_7,l_8,l_9$, depending on $n$, such that $0\leq l_4, l_5,l_6 \leq p$ and
\begin{align}\label{eqn:tdonc0}
&t_1\leq t_2 \leq \dots \leq t_{l_4} < -\epsilon^{1/3} \omega_{n}\nonumber\\
-\epsilon^{1/3}\omega_{n} \leq\; & t_{l_4+1}\leq t_{l_4+2}\leq \dots \leq\; t_{l_5} < 1\nonumber\\
1 \leq\; & t_{l_5+1}\leq  t_{l_5+2}\leq \dots \leq t_{l_6} \leq \epsilon^{1/3} \omega_{n}\nonumber\\
\epsilon^{1/3}\omega_{n} <\; & t_{l_6+1}\leq  t_{l_6+2}\leq \dots \leq t_p .
\end{align}
Use the same definitions of $l(x)$ and $d(x)$, as in the proof of Theorem \ref{theo:1}. Note that using Lemma 1 of C\"org\H{o} et al. (2003), we have $\delta_n\rightarrow 0$ as $n\rightarrow \infty$. Then due to Lemma \ref{lem:1}, \ref{lem:5} and the fact that $\epsilon \leq c < 1/8$, we have
\begin{align}\label{eqn:tdonc2}
d(x) \leq\; \Big[C_3\delta_ne^{-x^2/2}\Big]I\Big(2|x|\leq 2\Big) + \Big[C_3\epsilon^{5/3} \omega_{n}^{-1}e^{-x^2/2}\Big]I\Big(1<|x|\leq \epsilon^{1/3}\omega_n\Big)
\leq\;\Big[C_3\epsilon^{5/3} \omega_{n}^{-1}e^{-x^2/2}\Big],
\end{align}
when $|x|\leq \epsilon^{1/3}\omega_n$, and
\begin{align}\label{eqn:tdonc22}
d(x) 
& \leq \mathbf{P}\Big(T_{n1} \leq x \Big) + \Phi(x)\nonumber\\
&  \leq \Big|F_n\big(-\epsilon^{1/3}\omega_n\big)-\Phi\big(-\epsilon^{1/3}\omega_n\big)\Big| + \Big|\big(1-F_n\big(\epsilon^{1/3}\omega_n\big)\big)-\big(1-\Phi\big(\epsilon^{1/3}\omega_n\big)\big)\Big|+ 4\big(1-\Phi\big(\epsilon^{1/3}\omega_n\big)\big)\nonumber\\
& \leq 2\Big[C_3\epsilon^{5/3} \omega_{n}^{-1}e^{-x^2/2}\Big] + 4\phi(\epsilon^{1/3}\omega_n)/(\epsilon^{1/3}\omega_n)\leq 2e^{-(\epsilon^{1/3}\omega_n^{2})/2},
\end{align}
when $|x|> \epsilon^{1/3}\omega_n$. Again using Lemma \ref{lem:1} and \ref{lem:5} for sufficiently large $n$ we have
\begin{align}\label{eqn:tdonc1}
l(x) \leq\; & I\Big(x>\epsilon^{1/3}\omega_n\Big) + \bigg[1-\dfrac{2\phi(1)}{\sqrt{5}+1}+ C_3\delta_n e^{-1/2}\bigg]I\Big(x < 1\Big)\nonumber\\
& +\bigg[1-\dfrac{2\phi(x)}{\sqrt{x^2+4}+x}+C_3 \epsilon^{5/3}\delta_n^{1/6} e^{-x^2/2}\bigg]I\Big(x \in\Big[1,\epsilon^{1/3}\omega_n\Big]\Big)\nonumber\\
\leq\; & I\Big(x>\epsilon^{1/3}\omega_n\Big) + \bigg[1-\dfrac{\phi(1)}{\sqrt{5}+1}\bigg]I\Big(x < 1\Big)+\bigg[1-\dfrac{\omega_{n}^{-1}e^{-x^2/2}}{\sqrt{2\pi}\Big(\sqrt{\epsilon^{2/3}+4\omega_{n}^{-2}}+\epsilon^{1/3}\Big)} \bigg]I\Big(x \in\Big[1,\epsilon^{1/3}\omega_n\Big]\Big),
\end{align}
and
 Therefore from equations (\ref{eqn:tdonc0})-(\ref{eqn:tdonc22}) we have for sufficiently large $n$,
\begin{align}\label{eqn:tdonc3}
L(\bm{t})\leq J_5(\bm{t}) +J_6(\bm{t}) +J_7(\bm{t}) + J_8(\bm{t}),
\end{align}
where similar to the proof of Theorem \ref{theo:1}, we have
\begin{align}\label{eqn:tdonc4}
 J_5(\bm{t})+J_8(\bm{t}) \leq p\Big(2e^{-2^{-1}\epsilon^{2/3}\omega^2_n}\Big)   
 \leq 2\exp\Big(-\epsilon^{2/3}\omega_n^2(1/2-\epsilon^{1/3})\Big)
 < \epsilon^{5/3}/12,
\end{align}
for large enough $n$, since $\epsilon < 1/8$. Again similar to $J_2(\bm{t})$, it can be shown that 
\begin{align}\label{eqn:tdonc5}
 J_6(\bm{t})  \leq \Big(d^{-1} (\log d)^{-1}d^{-(\log d)^{-1}}\Big) C_3\epsilon^{5/3}\omega_n^{-1} < \epsilon^{5/3}/12, 
\end{align}
for sufficiently large $n$ where $d^{-1}=\Big[1-\frac{\phi(1)}{\sqrt{5}+1}\Big]$. The only thing that remains to bound is $J_7(\bm{t})$. Similar to the treatment to $J_3(\bm{t})$ in Theorem \ref{theo:1}, it can be shown that
\begin{align}\label{eqn:tdonc6}
J_7(\bm{t})\leq \epsilon^{5/3}/4 +\epsilon^{5/3}/12.
\end{align}
Therefore combining (\ref{eqn:tdonc3})-(\ref{eqn:tdonc6}), the proof is complete.\\

\noindent
{\bf Proof of Theorem \ref{theo:4}:} Here also we will follow the same route as in case of the proof of Theorem \ref{theo:1} or Theorem \ref{theo:3}, but obviously with different estimates of $l(x)$ and $d(x)$. Now take $c=\min\Big\{\big[12C_2\sqrt{2\pi}(\sqrt{2}+1)\big]^{-3},1/8\Big\}$ where $C_2$ is the constant defined in Lemma \ref{lem:4} with $Z_i$ replaced by $X_{i1}$. 
By Lemma \ref{lem:1}, Lemma \ref{lem:4} and the fact that $\epsilon \leq c < 1/8$, we have for sufficiently large $n$,
\begin{align}\label{eqn:tb1}
l(x) \leq\; & I\Big(x>\epsilon^{1/3}n^{1/4}\Big) + \bigg[1-\dfrac{2\phi(1)}{\sqrt{5}+1}+ 2C_2n^{-1/2} e^{-1/2}\bigg]I\Big(x < 1\Big)\nonumber\\
& +\bigg[1-\dfrac{2\phi(x)}{\sqrt{x^2+4}+x}+C_2\big(1+\epsilon^{1/3}n^{1/4}\big)n^{-1/2} e^{-x^2/2}\bigg]I\Big(x \in\Big[1,\epsilon^{1/3}n^{1/4}\Big]\Big)\nonumber\\
\leq\; & I\Big(x>\epsilon^{1/3}n^{1/4}\Big) + \bigg[1-\dfrac{\phi(1)}{\sqrt{5}+1}\bigg]I\Big(x < 1\Big)\nonumber\\
& +\bigg[1-\dfrac{n^{-1/4}e^{-x^2/2}}{\sqrt{2\pi}\Big(\sqrt{\epsilon^{2/3}+4n^{-1/2}}+\epsilon^{1/3}\Big)} \bigg]I\Big(x \in\Big[1,\epsilon^{1/3}n^{1/4}\Big]\Big),
\end{align}
and
\begin{align}\label{eqn:tb2}
d(x) \leq \Big[2 e^{-\big(\epsilon^{2/3}n^{1/2}\big)/2}\Big]I\Big(|x|>\epsilon^{1/3}n^{1/4}\Big)  + \Big[2C_2\epsilon^{1/3} n^{-1/4}e^{-x^2/2}\Big]I\Big(1<|x|\leq \epsilon^{1/3}n^{1/4}\Big).
\end{align}
Clearly (\ref{eqn:tb1}) and (\ref{eqn:tb2}) are respectively same as (\ref{eqn:tdonc1}) and (\ref{eqn:tdonc2}), but after replacing $\omega_n$ by $n^{1/4}$ and $C_3\epsilon^{5/3}$ by $3C_2\epsilon$. Therefore, the rest of the proof follows exactly following the arguments of the proof of Theorem \ref{theo:4}.\\

\noindent
{\bf Proof of Theorem \ref{theo:6}:}
We are going to prove the version of part (b) of Theorem \ref{theo:1} for $W_n$. Proof of the versions of Theorem \ref{theo:3} and \ref{theo:4} for $W_n$ are analogous. Although the proof will follow the steps as in the proof of part (b) of Theorem \ref{theo:1}, we need to take care the effect of the relation (\ref{eqn:relation}) on the estimates of $\tilde{l}(x)$ \& $\tilde{d}(x)$ (defined below) and then to utilize the monotonicity of the functions $f(y)=ny/(n+y-1)$ \& $g(y)=y\big[1-n(n+y-1)^{-1}\big]$ effectively to complete the proof. Note that we need to show
\begin{align}\label{eqn:tt1}
\limsup_{n\rightarrow \infty}\sup_{B\in\mathcal{A}^{re}}\Big|\mathbf{P}\Big(W_n\in B\Big)-\Phi\big(B\big)\Big|< \epsilon^{(1+\delta)/3},
\end{align}
whenever $\log p =\epsilon*d_{n,\delta}^2$ for some $\epsilon< c$ with some $0<c\leq 1/8$. Now note that the conclusion of lemma \ref{lem:7} is still true if we replace $D_n$ by $\tilde{D}_n=diag\Big(\tilde{d}_{n1},\dots,\tilde{d}_{np}\Big)$ in its statement where $\tilde{d}_{nj}=\frac{n}{n-1}\sum_{i=1}^{n}\Big(Z_{ij}-\bar{Z}_{nj}\Big)^2$ and $\bar{Z}_{nj}=n^{-1}\sum_{i=1}^{n}Z_{ij}$, $j\in \{1,\dots,p\}$. Therefore using Lemma \ref{lem:7}, to prove (\ref{eqn:tt1}) it is enough to show that
\begin{align}\label{eqn:tt2}
\sup_{t_1\leq t_2\leq \dots \leq t_p} \tilde{L}((t_1,\dots,t_p)^\prime)= \sup_{t_1\leq t_2\leq \dots \leq t_p}\bigg[\sum_{j=1}^{p}\Big(\prod_{j\neq k}\tilde{l}(t_j)\Big)\tilde{d}(t_k)\bigg]< \epsilon^{(1+\delta)/3}/2,
\end{align}
for large enough $n$. Here, $\tilde{l}(x)=\max\Big{\{}\mathbf{P}\Big(W_{n1} \leq x \Big), \mathbf{P}\Big(S_{n1} \leq x \Big)\Big{\}} \;\;\; \text{and}\;\;\; \tilde{d}(x)=\Big{|}\mathbf{P}\Big(W_{n1} \leq x \Big)-\mathbf{P}\Big(S_{n1} \leq x\Big)\Big{|}$, with $S_n=(S_{n1},\dots,S_{np})^\prime=n^{-1/2}\sum_{i=1}^{n}N_i$. Now fix $\bm{t}=(t_1,\dots,t_p)^\prime$ and consider the partition \ref{eqn:14} as in Theorem \ref{theo:1}. Next step is to find bounds on $\tilde{l}(x)$ and $\tilde{d}(x)$. Now from (\ref{eqn:relation}) we have $$\mathbf{P}\Big(W_{n1}\leq x\Big)=\mathbf{P}\Big(T_{n1}\leq x\Big(\frac{n}{n+x^2-1}\Big)^{1/2}\Big),$$ for any $x\in \mathcal{R}$. Therefore using Lemma \ref{lem:3}, for sufficiently large $n$ we have 
\begin{align}\label{eqn:tt3}
&\Big|\mathbf{P}\Big(W_{n1}\leq x\Big)-\Phi(x)\Big|\nonumber\\ \leq\;& \Big|\mathbf{P}\Big(T_{n1}\leq x\Big(\frac{n}{n+x^2-1}\Big)^{1/2}\Big)-\Phi\Big(x\Big(\frac{n}{n+x^2-1}\Big)^{1/2}\Big)\Big|+\Big|\Phi\Big(\Big(\frac{n}{n+x^2-1}\Big)^{1/2}\Big)-\Phi(x)\Big|\nonumber  \\
 \leq\;& A(1+2|x|)^{1+\delta}\exp\Big(-2^{-1}x^2\Big(\frac{n}{n+x^2-1}\Big)\Big)d_{n,\delta}^{-(2+\delta)}\nonumber\\
 &+ 3(\sqrt{2\pi})^{-1}\Big(|x|^3+|x|\Big)\exp\Big(-2^{-1}x^2\Big(\frac{n}{n+x^2-1}\Big)\Big)[n(n-1)]^{-1/2}\nonumber\\
\leq\;&  3A\big(1+2|x|\big)^{1+\delta}\exp\Big(-2^{-1}x^2\Big(\frac{n}{n+x^2-1}\Big)\Big)d_{n,\delta}^{-(2+\delta)}\nonumber\\
 \leq\; & A3^{2+\delta}\epsilon^{(1+\delta)/3} d_{n,\delta}^{-1}\exp\Big(-2^{-1} \epsilon^{2/3}d_{n,\delta}^2\Big(\frac{n}{n+d_{n,\delta}^2-1}\Big)\Big),
\end{align}
whenever $|x|\leq d_{n,\delta}$, since $d_{n,\delta}\leq n^{1/4}$ for any $0< \delta\leq 1$. Again whenever $|x|\geq \epsilon^{1/3}d_{n, \delta}$ using Lemma \ref{lem:1} and Lemma \ref{lem:3} we have for sufficiently large $n$, 
\begin{align}\label{eqn:tt4}
&\Big|\mathbf{P}\Big(W_{n1}\leq x\Big)-\Phi(x)\Big|\nonumber\\
& \leq 1-\mathbf{P}\Big(T_{n1}\leq \epsilon^{1/3}d_{n,\delta}\Big(\frac{n}{n+d_{n,\delta}^2-1}\Big)^{1/2}\Big) + \mathbf{P}\Big(T_{n1}\leq -\epsilon^{1/3}d_{n,\delta}\Big(\frac{n}{n+d_{n,\delta}^2-1}\Big)^{1/2}\Big)\nonumber\\
&\;\;\;+ 1-\Phi\Big(\epsilon^{1/3}d_{n,\delta}\Big(\frac{n}{n+d_{n,\delta}^2-1}\Big)^{1/2}\Big)
\leq 2\exp\Big(-2^{-1} \epsilon^{2/3}d_{n,\delta}^2\Big(\frac{n}{n+d_{n,\delta}^2-1}\Big)\Big).
\end{align}


Now define $8c=\min\big\{\big[A3^{2+\delta}\sqrt{2\pi}(\sqrt{2}+1)\big]^{-3},1\big\}$ and consider the partition of $\{t_1,\dots,t_p\}$ as in the proof of Theorem \ref{theo:1}. Write $x_1=x\Big(\frac{n}{n+x^2-1}\Big)^{1/2}$. Then due to (\ref{eqn:tt3}) we have for sufficiently large $n$,
\begin{align}\label{eqn:tt5}
\tilde{l}(x)  \leq\; & I\Big(x>\epsilon^{1/3}d_{n,\delta}\Big) + \bigg[1-\dfrac{2\phi(1)}{\sqrt{5}+1}+ A3^{2+\delta} d_{n,\delta}^{-(2+\delta)}\bigg]I\Big(x < 1\Big)\nonumber\\
& +\bigg[1-\dfrac{2\phi(x)}{\sqrt{x^2+4}+x}+A3^{2+\delta}\epsilon^{(1+\delta)/3}d_{n,\delta}^{-1} e^{-x_1^2/2}\bigg]I\Big(x \in\Big[1,\epsilon^{1/3}d_{n,\delta}\Big]\Big)\nonumber\\
\leq\; & I\Big(x>\epsilon^{1/3}n^{1/4}\Big) + \bigg[1-\dfrac{\phi(1)}{\sqrt{5}+1}\bigg]I\Big(x < 1\Big)+\bigg[1-\dfrac{(2\pi)^{-1/2}d_{n,\delta}^{-1}e^{-x_1^2/2}}{\Big(\sqrt{\epsilon^{2/3}+4d_{n,\delta}^{-2}}+\epsilon^{1/3}\Big)} \bigg]I\Big(x \in\Big[1,\epsilon^{1/3}d_{n,\delta}\Big]\Big).
\end{align}
The third part of the second inequality follows due to the facts that  $A3^{2+\delta}\sqrt{2\pi}(\sqrt{2}+1)\epsilon^{(2+\delta)/3}< 1/8$ \& $x_1\leq x$ for $x\geq 1$ and by noting that $g(y)=y\big[1-n(n+y-1)^{-1}\big]$ is non-decreasing when $y\geq 1$. Again due to (\ref{eqn:tt3}) and (\ref{eqn:tt4}) we have
\begin{align}\label{eqn:tt6}
\tilde{d}(x) &\leq \bigg[2\exp\bigg(- \dfrac{\epsilon^{2/3}d_{n,\delta}^2}{2}\bigg(\frac{n}{n+d_{n,\delta}^2-1}\bigg)\bigg)\bigg]I\Big(|x|>\epsilon^{1/3}d_{n,\delta}\Big)+ \Big[\dfrac{A3^{2+\delta}\epsilon^{(1+\delta)/3}}{d_{n,\delta}e^{x_1^2/2}} \Big]I\Big(|x|\leq \epsilon^{1/3}d_{n,\delta}\Big)
\end{align}
for any $x\in \mathcal{R}$, for sufficiently large $n$. Therefore from equations (\ref{eqn:tt5}) and (\ref{eqn:tt6}) we have for sufficiently large $n$,
\begin{align}\label{eqn:tt7}
\tilde{L}(\bm{t})\leq \tilde{J}_1(\bm{t}) +\tilde{J}_2(\bm{t}) +\tilde{J}_3(\bm{t}) + \tilde{J}_4(\bm{t}),
\end{align}
where $\{\tilde{J}_i\}_{i=1}^{4}$ are same as $\{J_i\}_{i=1}^{4}$ in the proof of Theorem \ref{theo:1}, but after replacing $e^{-\epsilon^{2/3}d_{n,\delta}^{2}/2}$ by $\exp\Big(-2^{-1} \epsilon^{2/3}d_{n,\delta}^2\Big(\frac{n}{n+d_{n,\delta}^2-1}\Big)\Big)$ and $t_j$ by $t_{j1}=t_j\Big(\frac{n}{n+t_j^2-1}\Big)^{1/2}$ for all $j\in \{1,\dots,p\}$.

Therefore $\tilde{J}_1(\bm{t})$, $\tilde{J}_2(\bm{t})$ and $\tilde{J}_4(\bm{t})$ can be dealt with exactly similarly as for $J_1(\bm{t})$, $J_2(\bm{t})$ and $J_4(\bm{t})$ in the the proof of Theorem \ref{theo:1} and noting that $d_{n,\delta}\leq n^{1/4}$ for all $0< \delta \leq 1$. Hence enough to show $\tilde{J}_3(\bm{t})\leq \epsilon^{(1+\delta)/3}/3$ for sufficiently large $n$. Note that
\begin{align} \label{eqn:tt8}
\tilde{J}_{3}(\bm{t}) &\leq \bigg(\sum_{k=l_2+1}^{l_3}A3^{2+\delta}\epsilon^{(1+\delta)/3}d_{n,\delta}^{-1}e^{-t_{k1}^2/2}\bigg(\prod_{ \substack{j=l_2+1\\ j \neq k}}^{l_3}\bigg[1-\dfrac{d_{n,\delta}^{-1}e^{-t_{j1}^2/2}}{\sqrt{2\pi}\Big(\sqrt{\epsilon^{2/3}+4d_{n,\delta}^{-2}}+\epsilon^{1/3}\Big)} \bigg]\bigg)\bigg)=\tilde{J}_{31}(\bm{t})\;\;\; \text{(say)}.
\end{align}
and
\begin{align*}
\dfrac{\partial \tilde{J}_{31}(\bm{t})}{\partial t_l }= \bigg[A3^{2+\delta}\epsilon^{1+\delta}d_{n,\delta}^{-1}\Big[\frac{n(n-1)t_l}{n+t_l^2-1}\Big]e^{-t_{l1}^2/2}\prod_{\substack{j=l_2+1\\ j\neq l}}^{l_3}\tilde{z}_j\bigg]\times
\bigg[\sum_{\substack{k=l_2+1\\k\neq l}}^{l_3}\bigg(\frac{\tilde{z}_k}{1-\tilde{z}_k}\bigg)-1\bigg],
\end{align*}
where $\tilde{z}_k=d_{n,\delta}^{-1}e^{-t_{k1}^2/2}\Big[\sqrt{2\pi}\Big(\sqrt{\epsilon^{2/3}+4d_{n,\delta}^{-2}}+\epsilon^{1/3}\Big)\Big]^{-1}$, $k=1,\dots,p$.
Hence for any $l=l_2+1,\dots,l_3$, $\dfrac{\partial \tilde{J}_{31}(\bm{t})}{\partial t_l }\gtreqless 0$ if and only if 
\begin{align}\label{eqn:tt9}
\sum_{\substack{j=l_2+1 \\ j\neq l}}^{l_3}\dfrac{\tilde{z}_j}{1-\tilde{z}_j}\gtreqless 1,
\end{align}
Now using the fact that $f(y)=ny/(n+y-1)$ is non-decreasing for $y>0$, we can claim that $1> \tilde{z}_{l_2+1} \geq \dots \geq \tilde{z}_{l_3} >0$ for sufficiently large $n$, since $1 \leq t_{l_2+1}\leq \dots \leq t_{l_3}$. Hence for sufficiently large $n$, $$\dfrac{\tilde{z}_{l_2+1}}{1-\tilde{z}_{l_2+1}} \geq \dots \geq \dfrac{\tilde{z}_{l_3}}{1-\tilde{z}_{l_3}},$$ due to the fact that $z/(1-z)$ is increasing for $z\in (0,1)$. Therefore from (\ref{eqn:tt9}) we can say that $\tilde{J}_{31}(\bm{t})$ is non-increasing in $\{t_{l_2+1},\dots,t_m\}$ and non-decreasing in $\{t_{m+1},\dots,t_{l_3}\}$ where $(m-l_2)$ is a non-negative integer not more than $(l_3-l_2)$. Hence we can follow the steps which leads to (\ref{eqn:20}) in the proof of Theorem \ref{theo:1} and conclude the desired bound on $\tilde{J}_3(\bm{t})$. \\

\textbf{Proof of Theorem \ref{theo:8}}: Throughout this proof assume that $$a_n^{-2} = \dfrac{2(2+\kappa)(\log p d_{n, \kappa})}{d_{n, \kappa}^2}.$$ Hence it is enough to prove that 
\begin{align}\label{eqn:8.1}
\tau_{n, \mathcal{A}^{re}}, \gamma_{n, \mathcal{A}^{re}}\leq A2^{10 + \kappa} a_n^{-(2+\kappa)},
\end{align}
where the constant $A$ is as defined in Lemma \ref{lem:3} with $Z_i=X_{i1}$, $i\in \{1,\dots, n\}$. Here we will only present the proof of (\ref{eqn:8.1}) for $\tau_{n, \mathcal{A}^{re}}$. The proof for $\gamma_{n, \mathcal{A}^{re}}$ is similar.  Assume that $a_n^{2+\kappa} \geq A2^{10+\kappa}$. Otherwise 
(\ref{eqn:8.1}) is trivially true. 


Recall from the proof of Theorem \ref{theo:1} that $T_n=\Big(\frac{\sum_{i=1}^{n}X_{i1}}{\sqrt{\sum_{i=1}^{n}X_{i1}^2}},\dots,\frac{\sum_{i=1}^{n}X_{ip}}{\sqrt{\sum_{i=1}^{n}X_{ip}^2}}\Big)^\prime$ and $S_n=n^{-1/2}\sum_{i=1}^{n}N_i$. Let $T_n=(T_{n1},\dots,T_{np})^\prime$ and $S_n= (S_{n1},\dots,S_{np})^\prime$. Hence due to Lemma \ref{lem:7} it is enough to show that
\begin{align}\label{eqn:8.3}
L((t_1,\dots,t_p)^\prime)= \sum_{j=1}^{p}\Big(\prod_{j\neq k}l(t_j)\Big)d(t_k)\leq A2^{9 + \kappa} a_n^{-(2+\kappa)},
\end{align}
for any $\bm{t}=(t_1,\dots, t_p)^\prime \in \mathcal{R}^p$ such that $t_1\leq t_2 \leq \dots \leq  t_p$. Here, $l(x)=\max\Big{\{}\mathbf{P}\Big(T_{n1} \leq x \Big), \mathbf{P}\Big(S_{n1} \leq x \Big)\Big{\}} \;\;\; \text{and}\;\;\; d(x)=\Big{|}\mathbf{P}\Big(T_{n1} \leq x \Big)-\mathbf{P}\Big(S_{n1} \leq x\Big)\Big{|}$.

First note that $a_n = o(d_{n,\kappa})$. Now fix $\bm{t}=(t_1,\dots,t_p)^\prime$ in $\mathcal{R}^p$ such that $t_1\leq t_2 \leq \dots\leq t_p$. Then there exist integers $l_{11},l_{12},l_{13}$, depending on $n$, such that $0\leq l_7, l_8,l_9 \leq p$ and
\begin{align}\label{eqn:8.4}
&t_1\leq t_2 \leq \dots \leq t_{l_7} < -a_n^{-1}d_{n,\kappa}\nonumber\\
-a_n^{-1}d_{n,\kappa} \leq\; & t_{l_7+1}\leq t_{l_7+2}\leq \dots \leq\; t_{l_8} < 1\nonumber\\
1 \leq\; & t_{l_8+1}\leq  t_{l_8+2}\leq \dots \leq t_{l_9} \leq a_n^{-1} d_{n,\kappa}\nonumber\\
a_n^{-1}d_{n,\kappa} <\; & t_{l_9+1}\leq  t_{l_9+2}\leq \dots \leq t_p 
\end{align}
Now due to Lemma \ref{lem:1} and Lemma \ref{lem:3} and $a_n^{2+\kappa} \geq A2^{10+\kappa}$, for any $|x|\leq a_n^{-1}d_{n,\kappa}$ we have
\begin{align}\label{eqn:8.5}
d(x)  = \Big{|}\mathbf{P}\Big(T_{n1} \leq x \Big)-\mathbf{P}\Big(S_{n1} \leq x\Big)\Big{|}
\leq A2^{1+\kappa}a_n^{-(1+\kappa)}d_{n,\kappa}^{-1}e^{-x^2/2},
\end{align}
and for $|x|>a_n^{-1}d_{n,\kappa}$ we have
\begin{align}\label{eqn:5.6}
d(x)  
 & \leq \Big|F_n\big(-a_n^{-1}d_{n, \delta}\big)-\Phi\big(-a_n^{-1}d_{n, \delta}\big)\Big| + \Big|\big(1-F_n\big(a_n^{-1}d_{n, \delta}\big)\big)-\big(1-\Phi\big(a_n^{-1}d_{n, \delta}\big)\big)\Big|\nonumber\\
& \leq A2^{4+\kappa}e^{-(a_n^{-2}d_{n,\kappa}^{2})/2}.
\end{align}
Again Lemma \ref{lem:1} and Lemma \ref{lem:3} imply that
\begin{align}\label{eqn:8.7}
l(x) 
\leq\; & I\Big(x>a_n^{-1}d_{n,\kappa}\Big) + \bigg[1-\dfrac{\phi(1)}{\sqrt{5}+1}\bigg]I\Big(x < 1\Big)+\bigg[1-(16\pi)^{-1/2}a_n d_{n,\kappa}^{-1}e^{-x^2/2}\bigg]I\Big(x \in\Big[1,a_n^{-1}d_{n,\kappa}\Big]\Big),
\end{align}
for any $x\in \mathcal{R}$. 
Therefore from equations (\ref{eqn:8.5})-(\ref{eqn:8.7}), we have
\begin{align*}
L(\bm{t})\leq I_1(\bm{t}) +I_2(\bm{t}) +I_3(\bm{t}) + I_4(\bm{t}),
\end{align*}
where

\begin{align*}
&I_1(\bm{t}) = \bigg(\Big[1-\dfrac{\phi(1)}{\sqrt{5}+1}\Big]^{l_8-1}\bigg)*\bigg(\prod_{j=l_8+1}^{l_9}\Big[1-\dfrac{a_n d_{n,\kappa}^{-1}e^{-t_j^2/2}}{(16\pi)^{1/2}}\Big]\bigg)*\bigg(\sum_{k=1}^{l_7}A2^{4+\kappa}e^{-(a_n^{-2}d_{n,\kappa}^{2})/2}\bigg)*I\Big(l_7\geq 1\Big),\\
&I_2(\bm{t})=\bigg(\Big[1-\dfrac{\phi(1)}{\sqrt{5}+1}\Big]^{l_8-1}\bigg)*\bigg(\prod_{j=l_8+1}^{l_9}\Big[1-\dfrac{a_n e^{-t_j^2/2}}{(16\pi)^{1/2}d_{n,\kappa}}\Big]\bigg)*\bigg(\sum_{k=l_7+1}^{l_8}\dfrac{A2^{1+\delta}e^{-t_k^2/2}}{a_n^{(1+\kappa)}d_{n,\kappa}}\bigg)*I\Big((l_8-l_7) \geq 1\Big),\\
&I_3(\bm{t})=\bigg(\Big[1-\dfrac{\phi(1)}{\sqrt{5}+1}\Big]^{l_8}\bigg)*\bigg(\sum_{k=l_8+1}^{l_9}\dfrac{A2^{1+\kappa}e^{-t_k^2/2}}{a_n^{(1+\kappa)}d_{n,\kappa}}\Big(\prod_{ \substack{j=l_8+1\\ j \neq k}}^{l_9}\Big[1-\dfrac{a_n e^{-t_j^2/2}}{(16\pi)^{1/2}d_{n,\kappa}}\Big]\Big)\bigg)*I\Big((l_9-l_8) \geq 1\Big),\\
&I_4(\bm{t})=\bigg(\Big[1-\dfrac{\phi(1)}{\sqrt{5}+1}\Big]^{l_8}\bigg)*\bigg(\prod_{j=l_8+1}^{l_9}\Big[1-\dfrac{a_n e^{-t_j^2/2}}{(16\pi)^{1/2}d_{n\delta}}\Big]\bigg)*\bigg(\sum_{k=l_9+1}^{p}A2^{4+\kappa}e^{-(a_n^{-2}d_{n,\kappa}^{2})/2}\bigg)I\Big((p-l_9) \geq 1\Big).
\end{align*}

\emph{Bound on $I_1(\bm{t})+I_4(\bm{t})$}:
Note that $a_n^{-2} = \dfrac{2(2+\kappa)(\log p d_{n, \kappa})}{d_{n, \kappa}^2}$. Then we have
\begin{align}\label{eqn:8.8}
I_1(\bm{t})+I_4(\bm{t}) \leq 2p\Big(A2^{4+\kappa}e^{-(a_n^{-2}d_{n,\kappa}^{2})/2}\Big)\leq A2^{5+\kappa}e^{[\log p - (2+\kappa)(\log p d_{n, \kappa})]}\leq  A2^{5+\kappa}a_n^{-(2+\kappa)}
\end{align}
\emph{Bound on $I_2(\bm{t})$}:
Let $d^{-1}=\Big[1-\dfrac{\phi(1)}{\sqrt{5}+1}\Big]$. Then
\begin{align}
I_2(\bm{t}) &\leq \bigg(\Big[1-\dfrac{\phi(1)}{\sqrt{5}+1}\Big]^{l_{8}-1}\bigg)\bigg(\sum_{k=l_7+1}^{l_{8}}A2^{1+\kappa}a_n^{-(1+\kappa)}d_{n,\kappa}^{-1}e^{-t_k^2/2}\bigg)I\Big((l_{8}-l_{7}) \geq 1\Big)\nonumber\\
& \leq A2^{1+\kappa} d a_n^{-(1+\kappa)}d_{n,\kappa}^{-1} \Big[\sup_{x>0}\big(xd^{-x}\big)\Big]\nonumber\\
&\leq \Big(A2^{1+\kappa} d (\log d)^{-1}d^{-(\log d)^{-1}}\Big) a_n^{-(2+\kappa)} \nonumber\\
& \leq A 2^{8+\kappa}a_n^{-(2+\kappa)}.
\end{align}
\emph{Bound on $I_3(\bm{t})$}: Note that if $(l_9-l_8) = 0$ then $I_3(\bm{t})=0$ and there is nothing more to do. Hence assume $(l
_9-l_8)\geq 1$. Then we have
\begin{align} \label{eqn:8.9}
I_3(\bm{t}) \leq \;& \bigg(\sum_{k=l_8+1}^{l_9}A2^{1+\kappa}a_n^{-(1+\kappa)}d_{n,\kappa}^{-1}e^{-t_k^2/2}\Big(\prod_{ \substack{j=l_8+1\\ j \neq k}}^{l_9}\Big[1-(16\pi)^{-1/2}a_n d_{n\kappa}^{-1}e^{-t_j^2/2}\Big]\Big)\bigg) = I_{31}(\bm{t})\;\;\; \text{(say)}.
\end{align}

Now due to the similar reasoning as in Theorem \ref{theo:1}, from (\ref{eqn:8.9}) we have
\begin{align*}
I_3(\bm{t})\leq I_{31}((\bm{t}^{(1)\prime},\bm{t}^{(2)\prime})^\prime)
\end{align*}
where $\bm{t}^{(1)}$ is an $(m-l_{8})\times 1$ vector with each component being $1$ and $\bm{t}^{(2)}$ is an $(l_9-m)\times 1$ vector with each component being $a_n^{-1}d_{n,\kappa}$. Therefore using $a_n^{-2} = \dfrac{2(2+\kappa)(\log p d_{n, \kappa})}{d_{n, \kappa}^2} \leq A^{-1}2^{9+\kappa}$, we have 
\begin{align}\label{eqn:8.10}
I_3(\bm{t}) \leq\; & (m-l_8)\Big[1-(16\pi)^{-1/2}a_n d_{n,\kappa}^{-1}e^{-1/2}\Big]^{m-l_8-1}\big(A2^{1+\kappa}a_n^{-(1+\kappa)}d_{n,\kappa}^{-1}e^{-1/2}\big)\nonumber\\
&+ (l_9-m)\big(A2^{1+\kappa}a_n^{-(1+\kappa)}d_{n,\kappa}^{-1}e^{-(2^{-1}a_n^{-2}d_{n,\kappa}^2)}\big)\nonumber\\
\leq \;& 2\Big[\exp\Big(\log (m-l_8)-(m-l_8)\big((16\pi)^{-1/2}a_n d_{n,\kappa}^{-1}e^{-1/2}\big)\Big)\Big]\big(A2^{1+\kappa}a_n^{-(1+\kappa)}d_{n,\kappa}^{-1}e^{-1/2}\big)\nonumber\\
& + p\big(A2^{1+\kappa}a_n^{-(1+\kappa)}d_{n,\kappa}^{-1}e^{-(2^{-1}a_n^{-2}d_{n,\kappa}^2)}\big)\nonumber\\
\leq\;&  2\Big[\exp\Big(\sup_{x>0}\big[\log x-x\big((16\pi)^{-1/2}a_n d_{n,\kappa}^{-1}e^{-1/2}\big)\big]\Big)\Big]\big(A2^{1+\kappa}a_n^{-(1+\kappa)}d_{n,\kappa}^{-1}e^{-1/2}\big)\nonumber\\
&+ \Big(A2^{1+\kappa}a_n^{-(1+\kappa)}d_{n,\kappa}^{-1}e^{[\log p - (2+\kappa)(\log p d_{n, \kappa})]}\Big)\nonumber\\
 \leq\; & \big((16\pi)^{1/2}a_n^{-1} d_{n,\kappa}e^{1/2}\big)\big(A2^{1+\kappa}a_n^{-(1+\kappa)}d_{n,\kappa}^{-1}e^{-1/2}\big) + A2^{1+\kappa} a_n^{-(2+\kappa)}\nonumber\\
 \leq\;& \Big[1+(16\pi)^{1/2}\Big]A2^{1+\kappa}a_n^{-(2+\kappa)}
\end{align}
Combining (\ref{eqn:8.8}) - (\ref{eqn:8.10}), we can conclude (\ref{eqn:8.3}) and hence the proof of Theorem \ref{theo:8} is now complete.\\

\textbf{Proof of Theorem \ref{theo:7}}. Throughout this proof assume that $a_n^{-2} = \dfrac{13(\log p \omega_n)}{\omega_n^2},$ where $\omega_n = \delta_{n,0}^{-1/6}$. Hence it is enough to show that 
\begin{align}\label{eqn:7.1}
\tau_{n, \mathcal{A}^{re}}, \gamma_{n, \mathcal{A}^{re}}\leq M_{11} a_n^{-6},
\end{align}
for some constant $M_{11}>0$. Here we will only present the proof of (\ref{eqn:7.1}) for $\tau_{n, \mathcal{A}^{re}}$. Since $X_{11}$ is in the domain of attraction of the normal distribution, $\delta_{n,0}\rightarrow 0$ (i.e. $\omega_n\rightarrow \infty$) as $n\rightarrow \infty$ due to Lemma 1 of C\"org\H{o} et al. (2003). Let $N_1$ be a natural number such that for all $n\geq N_1$,
\begin{align}\label{eqn:7.2}
C_3a_n^{-5}\omega_n^{-1}+\dfrac{2}{\sqrt{2\pi}}a_n\omega_n^{-1}\leq 1,C_3\omega_n^{-6}e^{-1/2}\leq \dfrac{\phi(1)}{\sqrt{5}+1},a_n \geq 1 \; \text{and}\;
\log \omega_n \geq 2,
\end{align}
where the constant $C_3$ is as defined in Lemma \ref{lem:3}, but with $Z_i=X_{i1}$, $i\in \{1,\dots, n\}$.
For $n< N_1$, again (\ref{eqn:8.1}) is trivial, since $M_{11}$ can be taken sufficiently large. Therefore as in the proof of Theorem \ref{theo:8}, it is enough to prove (\ref{eqn:8.1}) for $\tau_{n, \mathcal{A}^{re}}$ when $n\geq N_1$. For this we will mostly follow the proof of Theorem \ref{theo:8}. The differences will be in the estimates of $l(x)$, $d(x)$ and in the partition of $t_1, \dots, t_p$ after fixing $\bm{t}= (t_1,\dots, t_p)^\prime$. In the partition, $d_{n, \kappa}$ will be replaced by $\omega_n$. Now due to (\ref{eqn:7.2}), Lemma \ref{lem:1}, Lemma \ref{lem:5}  and the fact that $a_n = o(\omega_n)$, we have
\begin{align}\label{eqn:7.4}
d(x) 
\leq\;  \Big[C_3a_n^{-5} \omega_{n}^{-1}e^{-x^2/2}\Big]I\Big(|x|\leq a_n^{-1}\omega_n\Big)+ 2e^{-(a_n^{-2}\omega_n^{2})/2}I\Big(|x|> a_n^{-1}\omega_n\Big),
\end{align}
whenever $n\geq N_1$.
Again using Lemma \ref{lem:1} and \ref{lem:5} for $n\geq N_1$ we have
\begin{align}\label{eqn:7.5}
l(x) 
\leq\; & I\Big(x>a_n^{-1}\omega_n\Big) + \bigg[1-\dfrac{\phi(1)}{\sqrt{5}+1}\bigg]I\Big(x < 1\Big)+\bigg[1-(4\pi)^{-1/2}a_n\omega_{n}^{-1}e^{-x^2/2}\bigg]I\Big(x \in\Big[1,a_n^{-1}\omega_n\Big]\Big).
\end{align}
Clearly the bounds obtained in (\ref{eqn:7.4}) and (\ref{eqn:7.5}) are similar to the bounds obtained in (\ref{eqn:8.5})-(\ref{eqn:8.7}). Therefore the rest of the proof follows exactly same way as in case of Theorem \ref{theo:8}.\\

\textbf{Proof of Proposition \ref{prop:2.1}}. Note that the distribution of $X_{11}$, mentioned above, is defined in the proof of Proposition 1 of Chistyakov and G\"otze (2004). 
Since $\dfrac{\log p}{d_{n, 1}^2}\rightarrow \infty$ and $\beta_n/\sqrt{n}\rightarrow 0$ as $n\rightarrow \infty$, let $\log p = h_n d_{n, 1}^2$ where $h_n \rightarrow \infty$ and $h_n\beta_n = o\big(\sqrt{n}\big)$ as $n\rightarrow \infty$. Note that $EX_{11}^2 =1$ and $E|X_{11}|^3 = (1+A_2\theta\sqrt{b})5\beta_n/3$ for some $|\theta|\leq 1$ and absolute constant $A_2>0$. Therefore $n^{1/6}/\beta_n^{1/3}$ and $d_{n, 1}$ have the same order and $d_{n, 1}\rightarrow \infty$ as $n\rightarrow \infty$. Again $h_n = o\Big(\sqrt{n}/\beta_n\Big) = o \Big(n^{2/3}/\beta_n^{4/3}\Big)$ which implies $h_n = o(d_{n, 1}^4)$ i.e. $\sqrt{\log p} = o(d_{n, 1}^3)$.

Now we are going to mostly follow the proof of Proposition 1.1 of Fan and Koike (2021). As in that proof, here also enough to show that there exists a sequence of real numbers $\{x_n\}_{n\geq 1}$ such that
\begin{align}\label{eqn:2.1.1}
\liminf_{n\rightarrow \infty}\bigg|P\Big(\max_{1\leq j \leq p}T_{nj}\leq x_n\Big) - P\Big(\max_{1\leq j \leq p}Z_{nj}\leq x_n\Big)\bigg|>0.
\end{align}
Let $x_n$ be a sequence of real numbers such that $[\Phi(x_n)]^p = e^{-1}$, i.e.
\begin{align}\label{eqn:2.1.2}
P\Big(\max_{1\leq j \leq p}Z_{nj}\leq x_n\Big) = e^{-1}.
\end{align}
Then we have $x_n/{\sqrt{2\log p}} \rightarrow 1$ and $p\big(1-\Phi(x_n)\big) \rightarrow 1$ as $n\rightarrow \infty$ (cf. Proof of Proposition 2.1 in Koike (2019)). Again by applying Theorem 1 in Arratia et al. (1989) with $I=\{1,\dots, p\}$, $B_{\alpha}= \{\alpha\}$ and $X_{\alpha}= I\{T_{n\alpha}>x_n\}$ we have
\begin{align}\label{eqn:2.1.3}
\bigg|P\Big(\max_{1\leq j \leq p}T_{nj}\leq x_n\Big) - e^{-\lambda_n}\bigg|\leq p \Big[P\Big(T_{n1}> x_n\Big)\Big]^2,
\end{align}
where $\lambda_n = p \Big[P\Big(T_{n1}> x_n\Big)\Big]$.
Therefore,
\begin{align}\label{eqn:2.1.4}
\bigg|P\Big(\max_{1\leq j \leq p}T_{nj}\leq x_n\Big) - P\Big(\max_{1\leq j \leq p}Z_{nj}\leq x_n\Big)\bigg| \geq \big|e^{-\lambda_n}-e^{-1}\big| - p \Big[P\Big(T_{n1}> x_n\Big)\Big]^2    
\end{align}
Note that $x_n = O(\log p) = o(d_{n, 1}^3)$. Hence considering $b$ to be sufficiently small, due to equations (2.1) and (2.2) of Chistyakov and G\"otze (2004) and noting that $L_n^{-1} = d_{n, 1}^3$ we have for sufficiently large $n$,
\begin{align}\label{eqn:2.1.5}
\lambda_n =  \Big[p\big(1-\Phi(x_n)\big)\Big] \dfrac{P\Big(T_{n1}> x_n\Big)}{1-\Phi(x_n)}
& \leq 2 e^{-A_3d_{n,1}^{-3}x_n^3}\nonumber\\
 &\leq 2 e^{-A_4d_{n,1}^{-3}(\log p)^{3/2}}\nonumber\\
& \leq 2 e^{-A_4d_{n,1}^{-3}(h_n d_{n, 1}^2)^{3/2}} \nonumber\\
& =  2 e^{-A_4h_n^{3/2}} \rightarrow 0, \;\text{as}\; n\rightarrow \infty,
\end{align}
for some absolute positive constants $A_3, A_4$. In the first, second and third inequality we have respectively used $p\big(1-\Phi(x_n)\big)\rightarrow 1$, $x_n/\sqrt{2\log p}$ $\rightarrow 1$ and $\log p = h_n d_{n, 1}^2$. Again due to (\ref{eqn:2.1.5}) we have $p \Big[P\Big(T_{n1}> x_n\Big)\Big]^2 = p^{-1}\lambda_n^2 \rightarrow 0$ as $n\rightarrow \infty$. Hence for sufficienly small $b$, from equations (\ref{eqn:2.1.2}) - (\ref{eqn:2.1.5}) we can say that for sufficiently large $n$,
$$\bigg|P\Big(\max_{1\leq j \leq p}T_{nj}\leq x_n\Big) - P\Big(\max_{1\leq j \leq p}Z_{nj}\leq x_n\Big)\bigg| \geq 1/2\big(e^{-1/2}-e^{-1}\big).$$ 
Therefore equation (\ref{eqn:2.1.1}) follows and the Proposition \ref{prop:2.1} is proved.\\

\noindent
{\bf Proof of Proposition \ref{prop:2.2}:} 
Note that $\sum_{i=1}^{n}X_{i1}^2=\dots=\sum_{i=1}^{n}X_{ip}^2=n$. Hence the proof of Proposition \ref{prop:2.2} follows essentially through the same line of the proof of Theorem 3 of \cite{DL}. The only difference is in the choice of the function $f(\cdot)$ in the form of the set $A = \Big(-\infty, n^{1/4}f(n)\Big]^p$. Here $\{f(n)\}_{n\geq 1}$ has to be considered to be a sequence such that $n^{3/4}f(n)$ is an even integer and $2\big(\log_3 n\big)^{1/4} \leq f(n)\leq \dfrac{n^{1/4}}{1+\eta}$ with some $\eta \in (0, 1)$ satisfying $\bigg[\frac{1}{(1+\eta)^3}+\frac{1}{7(1+\eta)^5}-1\bigg]\geq 0$. 
Then choosing
$\sqrt{n}[f(n)]^2\approx
2[\log p -\log n]$ when $2\sqrt{n\log n}< \log p< \Big[\frac{3\log n}{4}+\frac{n}{2(1+\eta)^2}\Big]$, $f(n)=\frac{n^{1/4}}{1+\eta}$ when $\Big[\frac{3\log n}{4}+\frac{n}{2(1+\eta)^2}\Big]\leq \log p \leq \log \Big[\sqrt{\frac{\pi}{2}}(\sqrt{n+4}+\sqrt{n})e^{n/2}\Big]$  and $f(n)=n^{1/4}$ when $p> \sqrt{\frac{\pi}{2}}(\sqrt{n+4}+\sqrt{n})e^{n/2}$, all the steps of the proof of Theorem 3 of \cite{DL} go through and the theorem follows.\\

\textbf{Proof of Proposition \ref{prop:4.1}}. Note that it is enough to show that there exists a sequence of real numbers $\{y_n\}_{n\geq 1}$ such that
\begin{align}\label{eqn:4.1.1}
\liminf_{n\rightarrow \infty}\sqrt{\dfrac{n}{(\log p)^3}}\bigg|P\Big(\max_{1\leq j \leq p}T_{nj}\leq y_n\Big) - P\Big(\max_{1\leq j \leq p}Z_{nj}\leq y_n\Big)\bigg|>0.
\end{align}
Let $y_n = x_n$ where $\{x_n\}_{n\geq 1}$ be the sequence of real numbers defined in the proof of Proposition \ref{prop:2.1}, i.e. $[\Phi(y_n)]^p = e^{-1}$. We also have $y_n/{\sqrt{2\log p}} \rightarrow 1$ and $p\big(1-\Phi(y_n)\big) \rightarrow 1$ as $n\rightarrow \infty$. Again by applying Theorem 1 in Arratia et al. (1989) with $I=\{1,\dots, p\}$, $B_{\alpha}= \{\alpha\}$ and $X_{\alpha}= I\{T_{n\alpha}>y_n\}$ we have
\begin{align}\label{eqn:4.1.2}
&\bigg|P\Big(\max_{1\leq j \leq p}T_{nj}\leq y_n\Big) - P\Big(\max_{1\leq j \leq p}Z_{nj}\leq y_n\Big)\bigg|\nonumber\\
\geq\;& \bigg|e^{-\lambda_{1n}} - e^{-\lambda_{2n}}\bigg| - p \Big[P\Big(T_{n1}> x_n\Big)\Big]^2 - p\Big[1-\Phi(y_n)\Big]^2,
\end{align}
where $\lambda_{1n} = p \Big[P\Big(T_{n1}> y_n\Big)\Big]$ and $\lambda_{2n} = p\Big[1-\Phi(y_n)\Big]$. Since $\lambda_{2n}\rightarrow 1$,  $p\Big[1-\Phi(y_n)\Big]^2 = O(p^{-1})$ as $n\rightarrow \infty$. Again noting that $y_n = O(\sqrt{\log p}) = o\big(n^{1/6}\big)$ and $d_{n, 1} $ is of order $n^{1/6}$, we have from Lemma \ref{lem:3} that $\lambda_{1n}/{\lambda_{2n}} \rightarrow 1$, implying that $p \Big[P\Big(T_{n1}> x_n\Big)\Big]^2 = O(p^{-1})$ as $n \rightarrow \infty$. Therefore due to the assumption that $\sqrt{n} = o\big(p (\log p)^{3/2}\big)$ we have 
\begin{align}\label{eqn:4.1.3}
 \sqrt{\dfrac{n}{(\log p)^3}} \Big[p \big[P\Big(T_{n1}> x_n\Big)\big]^2 + p\big[1-\Phi(y_n)\big]^2\Big] \rightarrow 0, \text{as}\; n \rightarrow \infty.
\end{align}
Now since $y_n = o(n^{1/6})$, by Theorem 1.2 of Wang (2005) we have
\begin{align}\label{eqn:4.1.4}
\dfrac{\lambda_{1n}}{\lambda_{2n}}\leq \bigg[e^{-\Big(\dfrac{y_n^3EX_{11}^3}{3\sqrt{n}}\Big)}\bigg]\Big[e^{-A_5\mathcal{L}_{n, y_n}}\Big],   
\end{align}
for some absolute constant $A_5 > 0$. Here for any $x > 0$, $\mathcal{L}_{n, x}$ is defined as 
\begin{align*}
\mathcal{L}_{n, x}= \;&(1+x)n^{-1/2}E|X_{11}|^3 + (1+ x)^3n^{-1/2}E|X_{11}|^3I\big(|X_{11}|>\sqrt{n}/(1+x)\big)\\
&+ (1+x)^4n^{-1}EX_{11}^4\big(|X_{11}|\leq \sqrt{n}/(1+x)\big).
\end{align*}
Note that for any $0< \epsilon< 1$,  we have
\begin{align*}
(1+y_n)^4EX_{11}^4\big(|X_{11}|\leq\; &  \sqrt{n}/(1+y_n)\big) \leq \epsilon(1+y_n)^3n^{-1/2} E|X_{11}|^3\\
&+ (1+y_n)^3n^{-1/2} E|X_{11}|^3I\big(|X_{11}|>\epsilon\sqrt{n}/(1+y_n)\big).
\end{align*}
Again $y_n/{\sqrt{\log p}}\rightarrow 1$, $ y_n = o(n^{1/6})$ and the condition $E|X_{11}|^3< \infty$ implies that $(1+y_n)n^{-1/2}E|X_{11}|^3 = o\big(n^{-1/2}y_n^3\big)$. Additionally the  dominated convergence theorem implies that $(1+ y_n)^3n^{-1/2}E|X_{11}|^3I\big(|X_{11}|>\sqrt{n}/(1+y_n)\big) = o\big(n^{-1/2}y_n^3\big)$. 
Therefore we have $\mathcal{L}_{n, y_n} = o\big(n^{-1/2}y_n^3\big)$. Hence by noting that $y_n = o(n^{1/6})$ and $EX_{11}^3 > 0$, from (\ref{eqn:4.1.4}) we have for large enough $n$,
\begin{align}\label{eqn:4.1.5}
&\dfrac{\lambda_{1n}}{\lambda_{2n}}\leq \Big(1 - \dfrac{y_n^3EX_{11}^3}{6\sqrt{n}}\Big)\Big(1 + o\Big(\dfrac{y_n^3}{\sqrt{n}}\Big)\Big)    \nonumber\\
\Rightarrow \;\;\;& \lambda_{1n} - \lambda_{2n} \leq -\dfrac{y_n^3EX_{11}^3}{12\sqrt{n}}\nonumber \\
\Rightarrow \;\;\;& e^{-(\lambda_{1n} - \lambda_{2n})} \geq  1 + \dfrac{y_n^3EX_{11}^3}{48\sqrt{n}}\nonumber\\
\Rightarrow \;\;\;& \sqrt{\dfrac{n}{(\log p)^3}}\Big|e^{-(\lambda_{1n} - \lambda_{2n})} - 1\Big| \geq  \dfrac{y_n^3EX_{11}^3}{48\sqrt{(\log p)^3}}\nonumber\\
\Rightarrow \;\;\;& \sqrt{\dfrac{n}{(\log p)^3}}\Big|e^{-\lambda_{1n}} - e^{-\lambda_{2n}}\Big| \geq  \dfrac{EX_{11}^3}{48 e} > 0.
\end{align} 
In the first line we have used the inequalities $e^{-x} \leq 1 -x/2$ and $e^x \leq 1 + \dfrac{7x}{4}$ both for $0\leq x \leq 1$. In the third line we have used the inequality $e^x \geq 1 + \dfrac{x}{4}$ for $0\leq x \leq 1$. And in the last line we have used the facts that $\lambda_{2n} \rightarrow 1$ and $y_n/\sqrt{2\log p} \rightarrow 1$ as $n\rightarrow \infty$. Now combining equations (\ref{eqn:4.1.2}), (\ref{eqn:4.1.3}) and (\ref{eqn:4.1.5}), the proof is complete.

\vskip 3mm
\noindent
{\bf Acknowledgements:}  The author would like to thank Prof. Shuva Gupta and Prof. S. N. Lahiri for many helpful discussions.

\end{document}